\documentclass[american,uppersorbian,english]{article}
\usepackage[T1]{fontenc}
\usepackage[latin2,latin9]{inputenc}
\usepackage{color}
\usepackage{babel}
\usepackage{array}
\usepackage{float}
\usepackage{url}
\usepackage{multirow}
\usepackage{amsmath}
\usepackage{amsthm}
\usepackage{amssymb}
\usepackage{graphicx}
\usepackage{esint}
\PassOptionsToPackage{normalem}{ulem}
\usepackage{ulem}
\usepackage[unicode=true,pdfusetitle,
 bookmarks=true,bookmarksnumbered=false,bookmarksopen=false,
 breaklinks=false,pdfborder={0 0 1},backref=false,colorlinks=false]
 {hyperref}
\usepackage{breakurl}

\makeatletter

\providecommand{\tabularnewline}{\\}
\floatstyle{ruled}
\newfloat{algorithm}{tbp}{loa}
\providecommand{\algorithmname}{Algorithm}
\floatname{algorithm}{\protect\algorithmname}

\newcommand{\lyxaddress}[1]{
\par {\raggedright #1
\vspace{1.4em}
\noindent\par}
}
\theoremstyle{plain}
\newtheorem{thm}{\protect\theoremname}

\@ifundefined{showcaptionsetup}{}{%
 \PassOptionsToPackage{caption=false}{subfig}}
\usepackage{subfig}
\makeatother

\addto\captionsamerican{\renewcommand{\algorithmname}{Algorithm}}
\addto\captionsamerican{\renewcommand{\theoremname}{Theorem}}
\addto\captionsenglish{\renewcommand{\theoremname}{Theorem}}
\addto\captionsuppersorbian{\renewcommand{\theoremname}{Theorem}}
\providecommand{\theoremname}{Theorem}

\begin{document}

\title{Efficient Computation of Null-Geodesic with Applications to Coherent Vortex Detection}

\author{Mattia Serra\thanks{Email: serram@ethz.ch} and George Haller\thanks{Corresponding author. Email: georgehaller@ethz.ch}}

\maketitle
\selectlanguage{american}%

\lyxaddress{\begin{center}
\textit{Institute for Mechanical Systems, ETH Zürich,\\ 
	Leonhardstrasse 21, 8092 Zurich, Switzerland}
\par\end{center}}

\noindent Recent results suggest that boundaries of coherent fluid vortices
(elliptic coherent structures) can be identified as closed null-geodesics
of appropriate Lorentzian metrics defined on the flow domain. Here we derive an automated method for computing such null-geodesics
based on the geometry of the underlying geodesic flow. Our approach simplifies and improves existing procedures for computing variationally defined Eulerian and Lagrangian vortex boundaries. As an illustration, we compute objective vortex boundaries from satellite-inferred
ocean velocity data. A MATLAB implementation of our method is available as supplementary material. 

\selectlanguage{american}%

\section{Introduction}

Typical trajectories of general unsteady flows show complex paths,
yet, their phase space often contains regions of organized behavior. Recent mathematical results offer a rigorous definition of Objective
Coherent Structures (OCSs), uncovering the skeletons behind these
well-organized regions. 

OCSs can be classified as Lagrangian Coherent Structures (LCSs) and Objective
Eulerian Coherent Structures (OECSs), depending on the time interval
over which they organize nearby trajectories. Specifically, LCSs \cite{LCSHallerAnnRev2015}
are influential over a finite time interval, while OECSs \cite{SerraHaller2015}
are infinitesimally short-term limits of LCSs. LCSs are suitable for understanding
and quantifying finite-time transport and mixing in fluid flows, intrinsically tied to a preselected time interval. OECSs, in contrast, can be computed at any time instant, and hence are free from any assumptions
on time scales. OECSs are, therefore, promising tools for flow-control and real-time decision-making problems \cite{Serra2016}. 

Among the different types of coherent structures, vortex-type (elliptic) structures
are perhaps the most relevant for transport prediction and estimation, as they carry the same fluid mass over extended distances.

Lagrangian coherent vortices, in the sense of \cite{BlackHoleHaller2013},
are encircled by elliptic LCSs, i.e., exceptional material barriers
that exhibit no appreciable stretching or folding over a finite
time interval. In contrast, Eulerian coherent
vortex boundaries (elliptic OECSs), in the sense of \cite{SerraHaller2015},
are the instantaneous limits of elliptic LCSs. As such, elliptic OECSs are distinguished curves characterized by a lack of short-term filamentation.
We will refer to elliptic LCSs and elliptic OECSs collectively as elliptic Objective Coherent Structures (OCSs). \textcolor{black}{An alternative method for the identification of Lagrangian coherent vortices can be found in \cite{Hadjighasem2016}. This method computes vortex boundaries as stationary curves of the underlying stretching-based variational problem numerically, as opposed to \cite{BlackHoleHaller2013} in which the variational problem is solved exactly.}

Variational arguments show that elliptic OCSs can be located as null-geodesics of suitably defined Lorentzian
metrics (\cite{BlackHoleHaller2013,SerraHaller2015}). Their computation, however, requires a number of non-standard steps that complicate its implementation. These steps include (i) an accurate computation of eigenvalues and eigenvectors
of tensor fields \cite{Farazmand2012a}; (ii) trajectory integration for direction
fields as opposed to vector fields \cite{Tchon2006}; (iii) detection of singularities (regions
of repeated eigenvalues) of tensor fields and identification of their topological type \cite{ShearlessBarrierFarazmand2014}; (iv) selection of Poincaré sections for locating closed direction-field trajectories (cf. \cite{AutoMdetectFlorian} or Appendix \ref{sec:DirFieldApproachIC}).  

We develop here a simple and accurate method for the computation
of closed null-geodesics in two dimensions as periodic solutions of the
initial value problem 

\begin{equation}
\label{eq:IVPFinal}
\begin{cases}
r^{\prime} = \mathcal{F}(r,A(r),\nabla A(r)),\ \ r:=\begin{bmatrix} x\\ \phi \end{bmatrix}\in\mathbb{R}^{2}\times\mathbb{S}^{1}, \\
\textcolor{black}{r(0) = r_{0},}
\end{cases}
\end{equation}
where $\mathcal{F}(r,A,\nabla A)$ denotes a three-dimensional vector field, and
$A\in\mathbb{R}^{2\times2}$ the metric tensor associated with the particular type of elliptic
OCSs. Based on topological properties of planar closed curves,
we also derive the set of admissible initial conditions $r_{0}$ for null-geodesics.
Seeking periodic orbits of the initial value problem \eqref{eq:IVPFinal} is a significant simplification over previous approaches that were designed to locate closed null-geodesics as closed orbits of non-orientable direction fields with a large number of singularities (see e.g, \cite{BlackHoleHaller2013,AutoMdetectFlorian} or Appendix \ref{sec:DirFieldApproach}). Specifically, Karrash et al. \cite{AutoMdetectFlorian} devised an automated scheme for the detection of null-geodesics which relies on locating tensor-field singularities \cite{DELMARCELLEPhDThesis1994}. The detection of such singularities, however, is a sensitive process. This sensitivity increases with the integration time, leading to artificial clusters of singularities (cf. \cite{AutoMdetectFlorian} or Fig. \ref{fig:ClustCGsing}), which in turns precludes the detection of the outermost coherent vortex boundaries. Our method overcomes these limitations and identifies closed null-geodesics of a general Lorentzian metric without restrictions on their geometry, or on the number and type of singularities present in their interior.

The global orientability of $\mathcal{F}(r,A,\nabla A)$ also allows for cubic
or spline interpolation schemes. This leads to more accurate computations
compared with the integration of direction fields, for which the lack
of global orientability necessitates the use of linear interpolation. These simplifications enable a fully automated and accurate detection of variationally defined vortex boundaries in any two-dimensional unsteady velocity field without reliance on user input.
The integration of the three-dimensional vector field \eqref{eq:IVPFinal},
as well as the computation of the admissible set of initial conditions
$r_{0}$, uses standard built-in MATLAB functions available as supplementary material to this paper. Finally, the
ODE in \eqref{eq:IVPFinal} can be used to compute null-geodesics of general Lorentzian metrics, and hence is
also relevant for hyperbolic and parabolic OCSs defined from variational principles in \cite{ShearlessBarrierFarazmand2014} and \cite{SerraHaller2015}. 


\section{Formulation of the problem\label{sec:FormulationPbl}}

We consider a variational problem

\begin{equation}
\begin{array}{cc}
Q[\gamma(s)]=\int_{\gamma}L(x(s),x^{\prime}(s))\,ds,\ \ \ \delta\int_{\gamma}L(x(s),x^{\prime}(s))\,ds=0,\end{array}\label{eq:Functional}
\end{equation}
with a quadratic Lagrangian
\begin{equation}
L(x,x^{\prime})=\frac{1}{2}\left\langle x^{\prime},A(x)x^{\prime}\right\rangle ,\label{eq:Lagrangian}
\end{equation}
where $A(x)$ is a tensor for all $x\in U\subset\mathbb{R}^{2}$,
and $\langle\cdot,\cdot\rangle$ denotes the Euclidean inner product.
Let $x:s\mapsto x(s),\ s\in[0,\sigma]\subset\mathbb{R}$,
denote the parametrization of a geodesic $\gamma$ of the metric 
\begin{equation}
\label{eq:Metric}
g_x(x^\prime,x^\prime)=\frac{1}{2}\left\langle x^{\prime},A(x)x^{\prime}\right\rangle,
\end{equation}
and $x'(s):=\tfrac{dx}{ds}$ its local tangent vector. 

The Euler--Lagrange equations \cite{Gelfand&FomCalcVar2000}
associated with \eqref{eq:Functional} are 
\[
\frac{1}{2}\nabla_{x}\left\langle x^{\prime},A(x)x^{\prime}\right\rangle -\frac{d}{ds}\left[A(x)x^{\prime}\right]=0,
\]

\noindent with the equivalent four-dimensional first-order formulation

\begin{equation}
\label{eq:LAgeq}
	\begin{aligned} 
x^{\prime} = & v,\\
v^{\prime} = & \frac{1}{2}A^{-1}(x)[\nabla_{x}\left\langle x^{\prime},A(x)x^{\prime}\right\rangle]-A^{-1}(x)[(\nabla_{x}A(x)v)v].
	\end{aligned} 
\end{equation}

\noindent Here, in tensor notation and with summation implied over
repeated indices,

\noindent 
\[
v^{\prime}_{i}=\frac{1}{2}A_{ij}^{-1}(x)v_{k}A_{kl,j}(x)v_{l}-A_{ij}^{-1}A_{jk,l}(x)v_{l}v_{k},\ i,j,k,l\in\text{\{}1,2\text{\}}.
\]

The functional $L(x,x^{\prime})$ in \eqref{eq:Lagrangian} has no explicit dependence on
the parameter $s$. By Noether's theorem \cite{Gelfand&FomCalcVar2000}, the metric $g_x(v,v)$ is a first integral for \eqref{eq:LAgeq}, i.e., 
\begin{equation}
g_{x(s)}(v(s),v(s)) =\frac{1}{2}\left\langle v(s),A(x(s))v(s)\right\rangle=g_0=\text{const.}.\label{eq:Noether}
\end{equation}
Any nondegenerate level surface satisfying $g_{x}(v,v)=g_{0}$
defines a three-dimensional invariant manifold for \eqref{eq:LAgeq}
in the four-dimensional space coordinatized by $(x,v)$. Differentiation with respect to $s$ along trajectories in this manifold gives

\begin{equation*}
\frac{dg_x}{ds}=\langle\nabla_{(x,v)}g_{x(s)}(v(s),v(s)),(x^{\prime},v^{\prime})\rangle = 0,
\end{equation*}
which is equivalent to
\begin{equation}
\label{eq:EulLagr_On_0LS-2}
\begin{aligned}
2\left\langle v^{\prime},A(x)v\right\rangle  & =-\left\langle v,(\nabla_{x}A(x)v)v\right\rangle,\\
 & =-v_{i}A_{ij,k}(x)v_{k}v_{j},\ i,j,k\in\text{\{}1,2\text{\}}.
\end{aligned}
\end{equation}
We denote with $(\cdot)_{\parallel}$ and $(\cdot)_{\perp}$ the components
of $(\cdot)$ along $v$ and $v^{\perp}=Rv$ respectively, where $R$ is a counterclockwise ninety-degree rotation matrix.
\noindent Expressing $v^{\prime}=v_{\parallel}^{\prime}+v_{\perp}^{\prime}$, we rewrite
equation \eqref{eq:EulLagr_On_0LS-2} as:

\begin{equation}
2\left\langle v_{\parallel}^{\prime},A(x)v\right\rangle +2\left\langle v_{\perp}^{\prime},A(x)v\right\rangle=-\left\langle v,(\nabla_{x}A(x)v)v\right\rangle .\label{eq:EulLagr_On_0LS-2-2}
\end{equation}
Of particular interest for us are null-geodesics of $g_x(v,v)$. Such curves satisfy $g_x(v,v)\equiv 0.$ In this case,
eq. \eqref{eq:EulLagr_On_0LS-2-2} simplifies to

\begin{equation}
2\left\langle v_{\perp}^{\prime},A(x)v\right\rangle=-\left\langle v,(\nabla_{x}A(x)v)v\right\rangle .\label{eq:EulLagr_On_0LS-2-2-1}
\end{equation}
This relationship holds in any dimension ($x\in\mathbb{R}^n$), but we keep our discussion two-dimensional to focus on coherent-structure detections in planar flows.

\section{Reduced three-dimensional null-geodesic flow \foreignlanguage{american}{\label{sec:FormulationPbln2}}}

\subsection{Flow reduction \label{sub:FlowReduction}}

We introduce polar coordinates
in the $v$ direction by letting
\begin{equation}
v=\rho e_{\phi},\qquad\rho\in\mathbb{R}^{+},\qquad e_{\phi}=\left(\cos\phi,\sin\phi\right)^{\top},\qquad\phi\in\mathbb{S}^{1},\label{eq:PolCoord}
\end{equation}
and rewrite eq. \eqref{eq:Noether} as

\begin{equation}
g_x(\rho e_\phi,\rho e_\phi)=\rho^2g_x(e_\phi, e_\phi)=g_{0}\equiv0\iff\frac{1}{2}\left\langle e_{\phi},A(x)e_{\phi}\right\rangle =0,\ \ \  x\in U,\ \phi\in\mathbb{S}^{1}.\label{eq:Noether-2-1}
\end{equation}
We also define the zero surface of $g_x$ as 

\begin{equation}
\mathcal{M}=\left\{ \left(x,\phi\right)\in U\times\mathbb{S}^{1}:\,\,g_x(e_\phi,e_\phi)=\frac{1}{2}\left\langle e_{\phi},A(x)e_{\phi}\right\rangle =0\right\}.
\label{eq:ZeroSetOfMetric}
\end{equation}
In addition, we rewrite eq. \eqref{eq:EulLagr_On_0LS-2-2-1}
as 

\begin{equation}
2\phi'\left\langle Re_{\phi},A(x)e_{\phi}\right\rangle = -\rho\left\langle e_{\phi},(\nabla_{x}A(x)e_{\phi})e_{\phi}\right\rangle ,\qquad R:=\left(\begin{array}{cc}
0 & -1\\
1 & 0
\end{array}\right),\label{eq:DerivingNoetherPC}
\end{equation}
or equivalently,

\begin{equation}
\begin{aligned}
x^{\prime} &=\rho e_{\phi},\\
\phi'& =-\rho\frac{\left\langle e_{\phi},(\nabla_{x}A(x)e_{\phi})e_{\phi}\right\rangle }{2\left\langle e_{\phi},R^{\top}A(x)e_{\phi}\right\rangle }.
\end{aligned}
\label{eq:3DreducedFlowAppoggio}
\end{equation}

Next, we rescale time along each trajectory $\left(x(s),\rho(s),\phi(s)\right)$ of \eqref{eq:3DreducedFlowAppoggio}
by letting
\begin{equation}
\bar{s}=\int_{0}^{s}\rho(\sigma)d\sigma,\label{eq:rescale-1-2-1}
\end{equation}
which gives $\frac{dx}{d\bar{s}}=e_{\phi},\ \ \frac{d\phi}{d\bar{s}}=-\frac{\left\langle e_{\phi},(\nabla_{x}A(x)e_{\phi})e_{\phi}\right\rangle }{2\left\langle e_{\phi},R^{\top}A(x)e_{\phi}\right\rangle }.$
We then drop the bar on $s$ to obtain the final form

\begin{equation}
\begin{aligned}
\frac{dx}{ds} & =e_{\phi},\\
\frac{d\phi}{ds} & =-\frac{\left\langle e_{\phi},(\nabla_{x}A(x)e_{\phi})e_{\phi}\right\rangle }{2\left\langle e_{\phi},R^{\top}A(x)e_{\phi}\right\rangle },
\end{aligned}
\label{eq:EulLagr_On_0LS_Pc-1}
\end{equation}
for the reduced three-dimensional null-geodesic flow, which is defined on the set
\begin{equation*}
V = \left\{ \left(x,\phi\right)\in U\times\mathbb{\mathbb{S}}^{1}:\ A(x)e_{\phi}\nparallel e_{\phi},\ A(x)\neq \mathbf{0}\right\} ,
\end{equation*}
where $\mathbf{0}\in\mathbb{R}^{2\times2}$ denotes the null tensor (cf. Appendix \ref{sec:DomExistGeodFlow}). In words, $V$ is the
set of points in $U\times\mathbb{S}^{1}$ where $A(x)$
is nondegenerate and $e_{\phi}$ is not aligned with the eigenvectors
of $A(x)$. Note that by construction, $\phi^{\prime}(s)$
is the pointwise curvature of $\gamma$.
An equation related to eq. \eqref{eq:EulLagr_On_0LS_Pc-1}
appears in \cite{Ying2006} for the geodesic flow associated with
the Riemannian metric on a general manifold, defined as the zero set of
a smooth function $F(x)$.

The ODE \eqref{eq:EulLagr_On_0LS_Pc-1}
has one dimension less than eq. \eqref{eq:LAgeq},
and the $x-$projection of its closed orbits coincide with closed null-geodesics on $(U,g_x)$. This follows from the equivalence of null-surfaces and null-geodesics in two dimensions.
In Appendix \ref{sec:HamiltReduction}, using the
Hamiltonian formalism, we derive an equivalent reduced geodesic flow
in the $(x,p)$ variables, with $p$ denoting the generalized momentum. 

\begin{figure}[h]
\subfloat[\selectlanguage{english}%
\selectlanguage{american}%
]{\includegraphics[width=0.3\columnwidth]{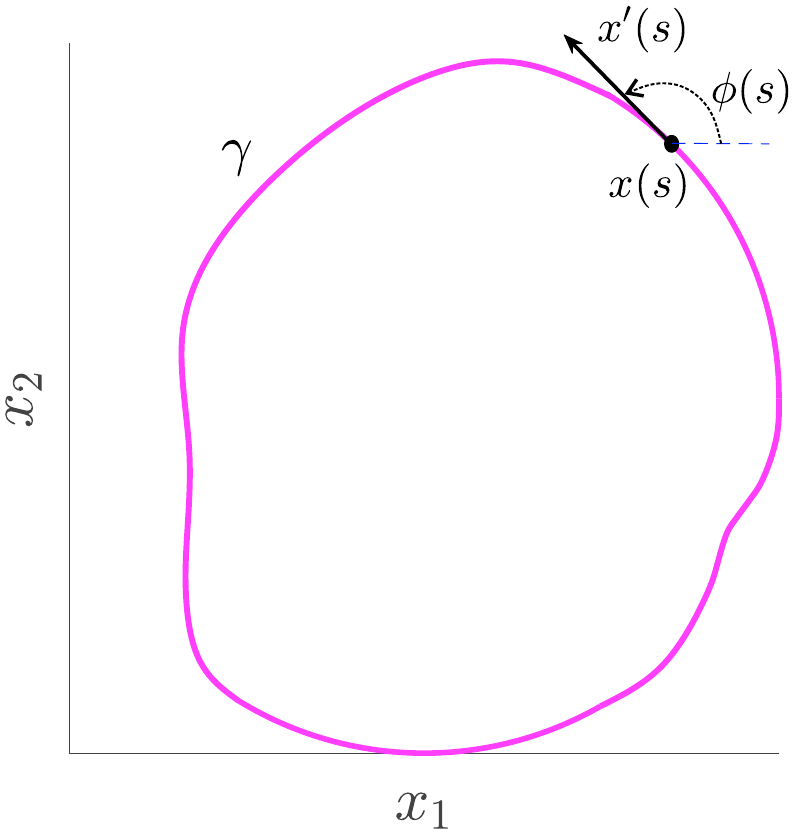}\label{fig:gamma2d}}\hfill{}\subfloat[\selectlanguage{english}%
\selectlanguage{american}%
]{\includegraphics[width=0.5\columnwidth]{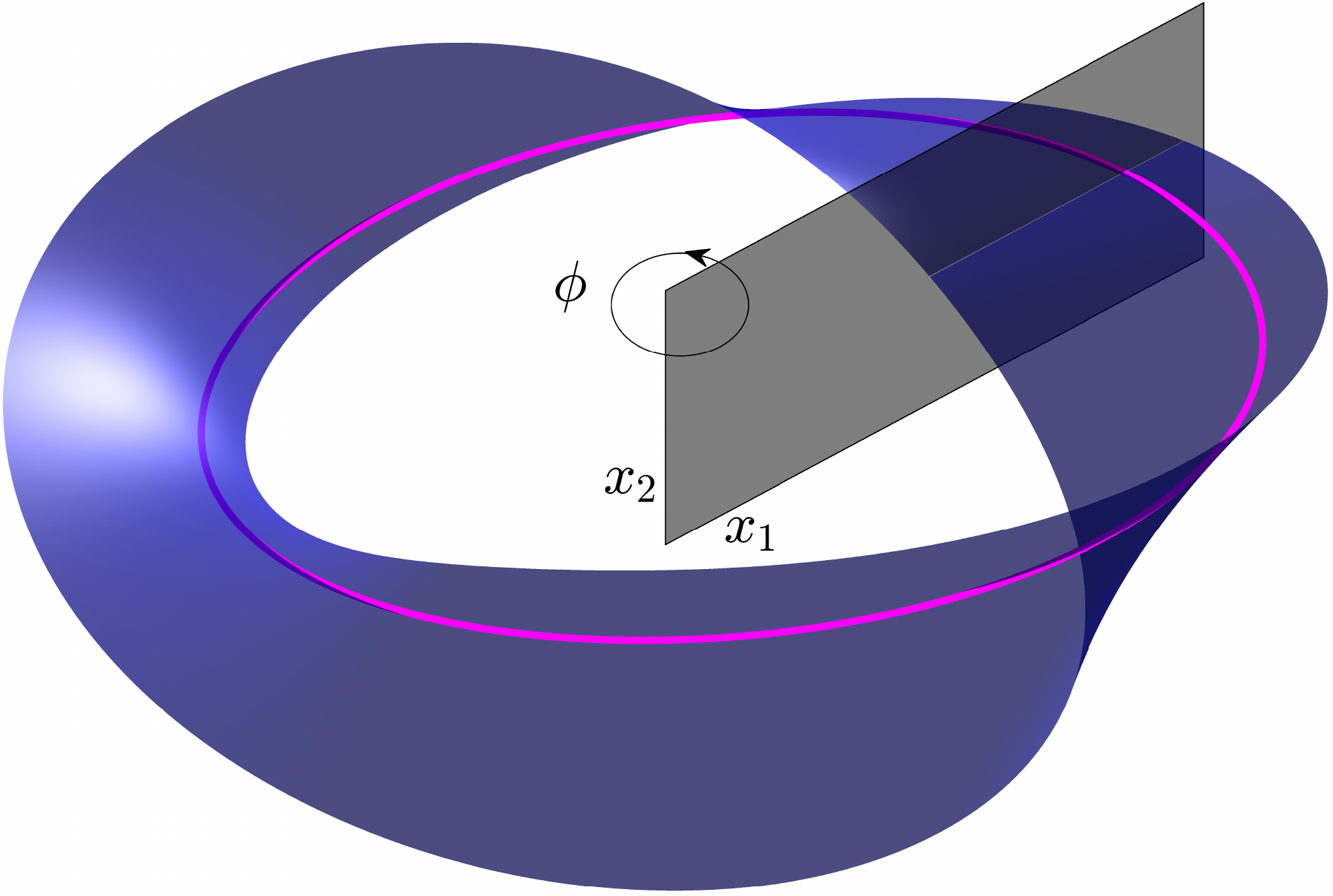}

\label{fig:gamma3d}}

\caption{(a) Closed null-geodesic of $g_x(u,u)$ in the $x-$subspace. (b) Closed
null-geodesic of $g_x(u,u)$ in $U\times\mathbb{\mathbb{S}}^{1}$ on
the zero level surface of $g_x(e_\phi,e_\phi)$. In this space, a
closed null-geodesics is an integral curve of \eqref{eq:EulLagr_On_0LS_Pc-1}
satisfying the boundary
conditions: $x(\sigma)=x(0),\ \text{and}\ \phi(\sigma)=\phi(0)\pm2\pi$.}

\label{fig:gamma2_3d}
\end{figure}

Figure \ref{fig:gamma2_3d} shows a closed null-geodesic
$\gamma$ of the metric $g_x(u,u)$ both in the $x-$subspace (Fig. \ref{fig:gamma2d})
and in the $U\times\mathbb{\mathbb{S}}^{1}-$space (Fig.
\ref{fig:gamma3d}). Specifically, Fig.
\ref{fig:gamma3d} shows a closed integral curve
of \eqref{eq:EulLagr_On_0LS_Pc-1} on the manifold $\mathbb{\mathcal{M}}$.

\subsection{Dependence on parameters\label{sub:Dependence-on-parameters}}

\selectlanguage{english}%
In applications to coherent vortex detection (cf. Section \textcolor{black}{\ref{sec:CoherentVtxBound}}),
the tensor field $A(x)$ depends on a parameter $\alpha\in\mathbb{R}$, leading to a specific 
Lorentzian metric family of the form 

\selectlanguage{american}%
\begin{equation}
g_{x,{\alpha}}=\frac{1}{2}\langle u,A_{\alpha}(x)u\rangle,\ \ \ \ A_{\alpha}=A(x)-\alpha I.\label{eq:TensorFamily}
\end{equation}
The zero level set of the metric family is defined by
$g_{x,\alpha}(e_\phi,e_\phi)=0.$ Interestingly, however, the reduced
ODE \eqref{eq:EulLagr_On_0LS_Pc-1}
remains independent of $\alpha$ because
\begin{equation}
\frac{\left\langle e_{\phi},(\nabla_{x}A_{\alpha}(x)e_{\phi})e_{\phi}\right\rangle }{2\left\langle e_{\phi},R^{\top}A_{\alpha}(x)e_{\phi}\right\rangle } = \frac{\left\langle e_{\phi},(\nabla_{x}A(x)e_{\phi})e_{\phi}\right\rangle }{2\left\langle e_{\phi},R^{\top}A(x)e_{\phi}\right\rangle }.
\end{equation}
We summarize this result in the following theorem. 
\begin{thm}
\label{sec:ThmIndepParam} The reduced three dimensional null-geodesic
flow of the Lorentzian metric family $g_{x,{\alpha}}(u,u)=\frac{1}{2}\left\langle u,A_{\alpha}(x)u\right\rangle ,\ A_{\alpha}(x)=A(x)-\alpha I,\ \ \alpha\in\mathbb{R}$,
is independent of $\alpha$ and satisfies the differential equation

\begin{equation*}
\begin{aligned}
x^{\prime} & = e_{\phi},\\
\phi^{\prime} & = -\frac{\left\langle e_{\phi},(\nabla_{x}A(x)e_{\phi})e_{\phi}\right\rangle }{2\left\langle e_{\phi},R^{\top}A(x)e_{\phi}\right\rangle },
\end{aligned}
\end{equation*}
defined on the set 
\begin{equation*}
V=\left\{ \left(x,\phi\right)\in U\times\mathbb{\mathbb{S}}^{1}:\ A(x)e_{\phi}\nparallel e_{\phi},\ A(x)\neq \mathbf{0}\right\}.
\end{equation*}
\end{thm}
The ODE \eqref{eq:EulLagr_On_0LS_Pc-1} is independent of $\alpha$, and hence
all null-geodesics of the metric family $g_{x,\alpha}$ can be integrated under the same vector field,
as opposed to available direction field formulations that depend
on $\alpha$ (cf. Appendix \ref{sec:DirFieldApproach}, eq. \eqref{eq:diffieldODE}).
This property of the ODE \eqref{eq:EulLagr_On_0LS_Pc-1} further simplifies
the computation of null-geodesics of $g_{x,{\alpha}}(u,u)$.

\subsection{Initial conditions \foreignlanguage{american}{\label{sec:Correct Initial conditions}}}

\selectlanguage{american}%
The only missing ingredient for computing null-geodesics of $g_{x,\alpha}$,
is a set of initial conditions for the reduced null-geodesic flow
\eqref{eq:EulLagr_On_0LS_Pc-1}.
Here we derive the set of initial conditions $r_{0}\subset V$,
such that any null-geodesic of $g_{x,\alpha}$, necessarily contains a point in $r_{0}$. According to Sections \ref{sec:FormulationPbl}-\ref{sub:Dependence-on-parameters},
for any fixed value of $\alpha$, null-geodesics of $g_{x,\alpha}$
must lie on the zero level surface of $g_{x,\alpha}(e_\phi,e_\phi)$, i.e., on
\[
\mathcal{M}_{\alpha}=\left\{ \left(x,\phi\right)\in U\times\mathbb{S}^{1}:\,\,g_{x,\alpha}(e_\phi,e_\phi)=0\right\} .
\]
Furthermore, for every closed planar curve
$\gamma$, the angle $\phi$ between its local tangent vector and an arbitrary fixed direction (cf. Fig. \ref{fig:gamma2d})
assumes all vales in the interval $[0,2\pi]$. This simple topological
property of closed regular planar curves allows us to define the
admissible set of initial conditions for \eqref{eq:EulLagr_On_0LS_Pc-1}
as follows.

\noindent For every fixed $\alpha$, we compute the set of
initial conditions $r_{\alpha}(0)$ as 
\begin{equation}
r_{\alpha}(0)=\left\{ \left(x_{0},\phi_{0}\right)\in V:\ g_{x_0,\alpha}(e_{\phi_{0}},e_{\phi_{0}})=0,\ \ \forall\phi_{0}\in\mathbb{S}^{1}\right\} .\label{eq:r0alpha}
\end{equation}

\begin{figure}[h]
\subfloat[\selectlanguage{english}%
\selectlanguage{american}%
]{\includegraphics[width=0.5\columnwidth]{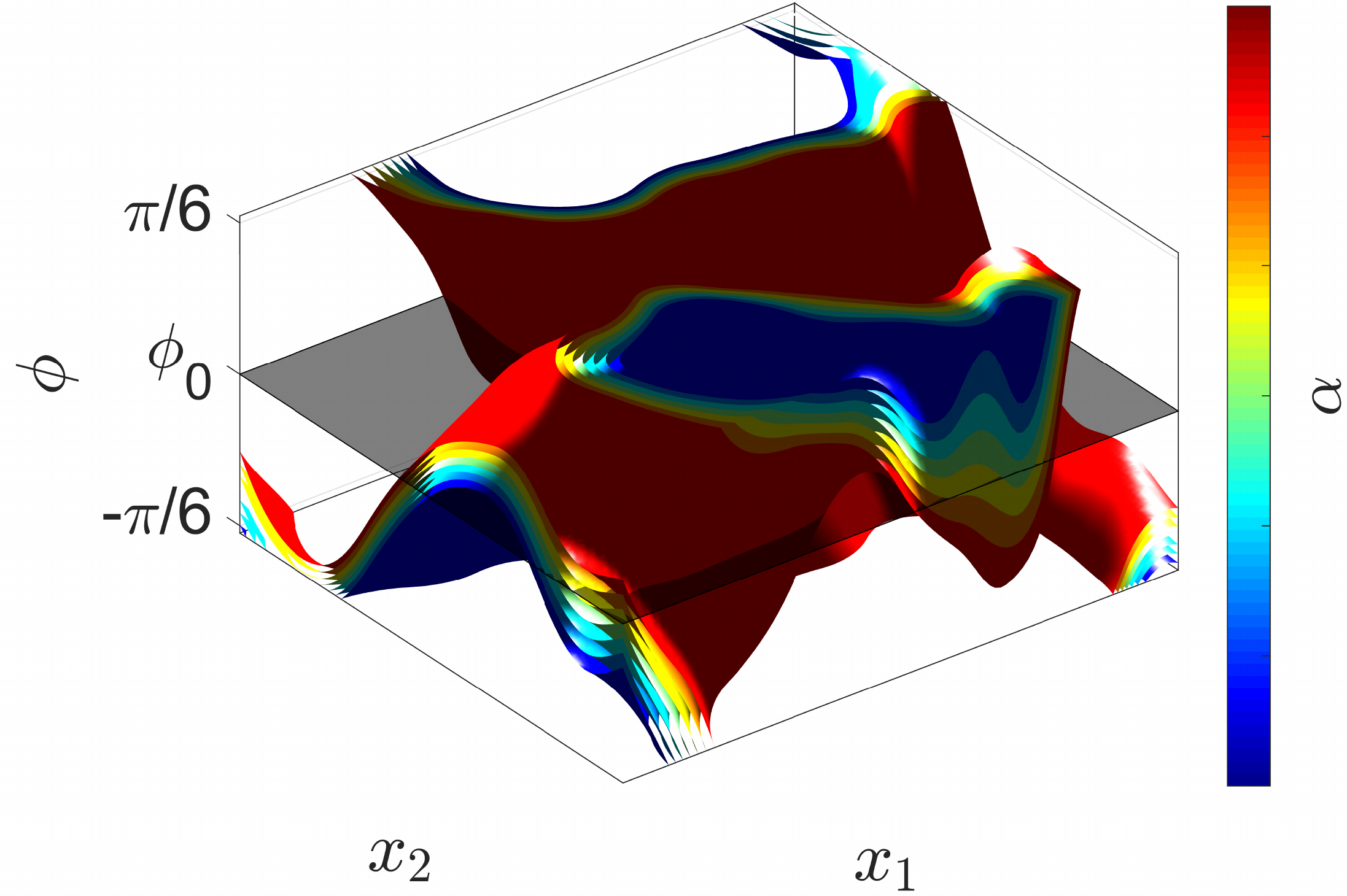}\label{fig:InitCondFoliat}}\hfill{}\subfloat[\selectlanguage{english}%
\selectlanguage{american}%
]{\includegraphics[width=0.4\columnwidth]{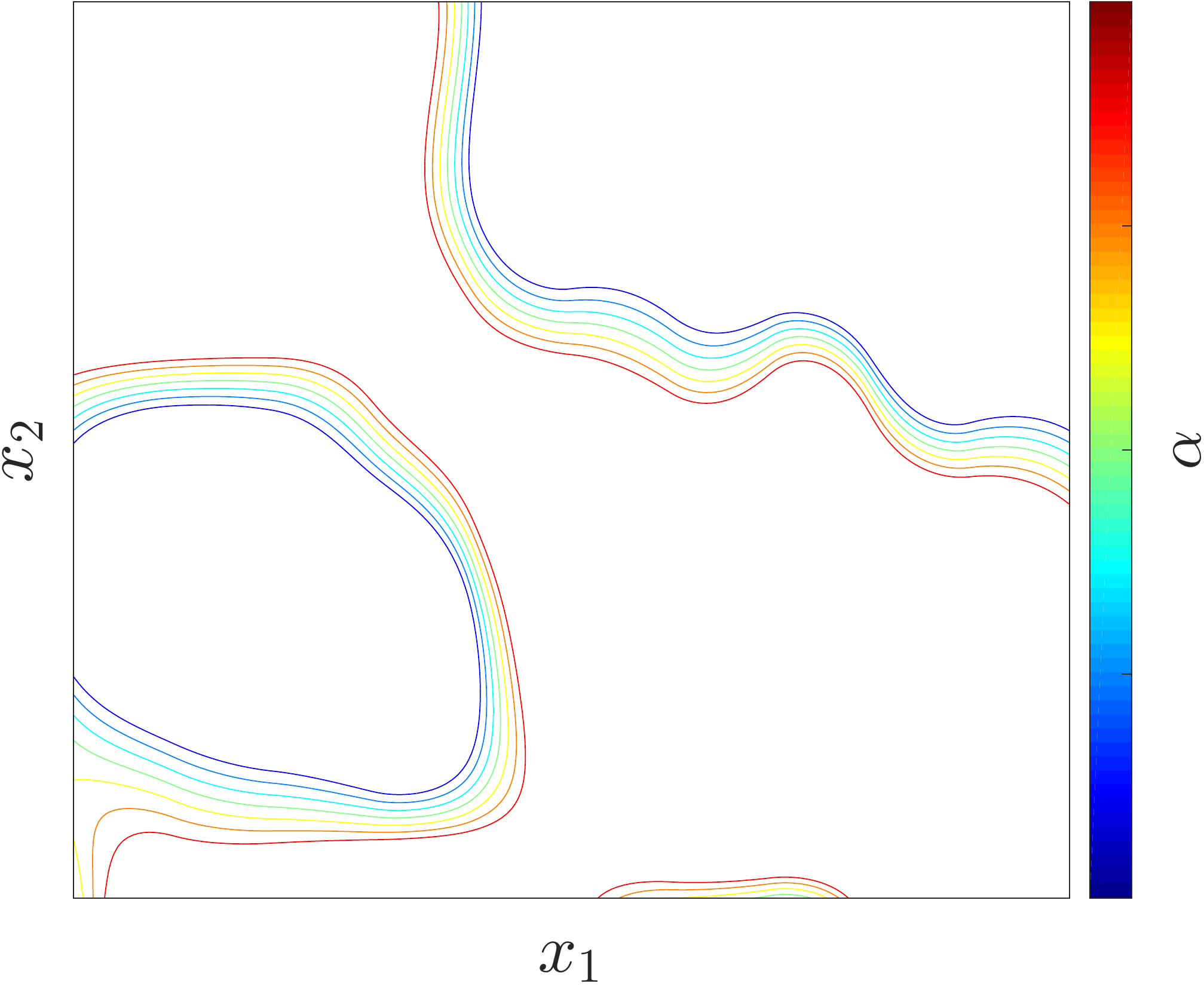}

\label{fig:InitCondCurves}}\caption{(a) Section of $\mathcal{M}_{\alpha}$ for different
values of $\alpha$, in
a flow example analyzed in more detail in Section \ref{sec:OceanExampleE}.
The black plane corresponds to $\phi(x)=\phi_{0}=0.$ (b) Set of points
$x_{0}(\alpha,\phi_{0})\in U$ satisfying $\mathcal{M}_{\alpha}\cap\phi(x)=\phi_{0}=0$,
for different values of $\alpha$.}
\label{fig:InitialCond}
\end{figure}
Figure \ref{fig:InitialCond}
illustrates formula \eqref{eq:r0alpha} in a flow example analyzed
in more detail in Section \ref{sec:OceanExampleE}. Specifically, Fig.
\ref{fig:InitCondFoliat} shows a section ($\phi\in[-\pi/6,\pi/6]$)
of $\mathcal{M}_{\alpha}$ for different values
of $\alpha$. The $\phi(x)=\phi_{0}=0$ plane is shown in black. Figure \ref{fig:InitCondCurves}
shows the set points $x_{0}(\alpha,0)\in U$ satisfying $g_{x_0,\alpha}(e_{0},e_{0})=0$.

\subsection{The initial value problem for null-geodesics }

Putting together the results from Sections \ref{sub:FlowReduction}-\ref{sec:Correct Initial conditions},
we obtain our main result, already summarized briefly in eq. \eqref{eq:IVPFinal}.
\begin{thm}
\label{sec:Thm1}Null-geodesics of the Lorentzian metric family $g_{x,\alpha}(u,u)=\frac{1}{2}\left\langle u,A_{\alpha}(x)u\right\rangle ,\ A_{\alpha}(x)=A(x)-\alpha I$,
coincide with the $x-$projection of closed orbits of the initial
value problem 
\begin{equation}
\begin{aligned}
x^{\prime} & = e_{\phi},\\
\phi^{\prime} & = -\frac{\left\langle e_{\phi},(\nabla_{x}A(x)e_{\phi})e_{\phi}\right\rangle }{2\left\langle e_{\phi},R^{\top}A(x)e_{\phi}\right\rangle },\label{eq:EulLagr_On_0LS_Pc-1-1}
\end{aligned}
\end{equation}
\begin{equation}
(x_{0},\phi_{0})=\left\{ \left(x_{0},\phi_{0}\right)\in V:\ \ \frac{1}{2}\langle e_{\phi_0},(A(x_0)-\alpha I)e_{\phi_0}\rangle=0,\ \ \forall\phi_{0}\in\mathbb{S}^{1}\right\} ,\label{eq:r0alpha-1}
\end{equation}
defined for any parameter value $\alpha\in\mathbb{R}$ on the set
\begin{equation}
V=\{x\in U,\ \phi\in\mathbb{S}^{1}:\ \left\langle e_{\phi},R^{\top}A(x)e_{\phi}\right\rangle \neq0\}.
\label{eq:DomDefRedFlow}
\end{equation}
Null-geodesics for a given value of $\alpha$, lie on the
zero level surface of $g_{x,\alpha}(e_\phi,e_\phi)$ defined as 
\begin{equation}
\mathcal{M}_{\alpha}=\left\{ \left(x,\phi\right)\in V:\ \frac{1}{2}\left\langle e_{\phi},(A(x)-\alpha I)e_{\phi}\right\rangle =0\right\} .\label{eq:LevelSurfaceThm}
\end{equation}
\end{thm}
\noindent \textcolor{black}{For any fixed valued of $\alpha$, the surface $\mathcal{M}_{\alpha}$ is a graph of the form $\phi(x,\alpha)$. Differentiating now $g_{x,\alpha}(e_{\phi(x,\alpha)},e_{\phi(x,\alpha))}$ with respect to $\alpha$, we obtain $\partial_\alpha\phi(x,\alpha)=\langle e_{\phi},R^{\top}A(x)e_{\phi}\rangle^{-1}$. This result, together with \eqref{eq:DomDefRedFlow} implies that null-geodesics of $g_{x,\alpha}$ cannot intersect for different values of $\alpha$, in agreement with the findings of \cite{BlackHoleHaller2013,SerraHaller2015}.}
\noindent In the following applications of Theorem \ref{sec:Thm1},
we select $\phi_{0}=0$ in formula \eqref{eq:r0alpha-1}.
 
\section{Null-Geodesics and the computation of Objective Coherent Structures}
\label{sec:CoherentVtxBound}

In the next section we recall the terminology used for the definition
of LCSs \cite{LCSHallerAnnRev2015} and OECSs \cite{SerraHaller2015}.

\subsection{Set-up and notation}\label{sec:setupAndNot}

Consider the two-dimensional non-autonomous dynamical system 
\begin{equation}
\dot{x}=f(x,t),\label{eq:FlowODE}
\end{equation}
with a twice continuously differentiable velocity field $f(x,t)$
defined over the open flow domain $U\in\mathbb{R}^{2}$, over a time
interval $t\in[a,b].$ We recall the customary velocity gradient decomposition
\[
\nabla f(x,t)=S(x,t)+W(x,t),
\]
with the rate-of-strain tensor $S=\tfrac{1}{2}(\nabla f+\nabla f^{\top})$
and the spin tensor $W=\tfrac{1}{2}(\nabla f-\nabla f^{\top})$. By our assumptions, both $S$ and $W$ are continuously differentiable in $x$ and $t$.

The rate-of-strain tensor is objective (i.e., frame-indifferent), whereas the spin tensor is
not objective as shown is classic texts on continuum mechanics (see,
e.g., \cite{TruesdellNoll2004}). The eigenvalues $s_{i}(x,t)$ and
eigenvectors $e_{i}(x,t)$ of $S(x,t)$ are defined, indexed and oriented
here through the relationship
\[
Se_{i}=s_{i}e_{i},\qquad\left|e_{i}\right|=1,\ \ i=1,2;\ \ s_{1}\leq s_{2},\quad e_{2}=Re_{1}.
\]
Fluid particles trajectories generated by $f(x,t)$ are solutions
of the differential equation $\dot{{x}}=f(x,t)$, and define
the flow map

\[
F_{t_{0}}^{t}(x_{0})=x(t;t_{0},x_{0}),\ \ \ \ x_{0}\in U,\ \ \ \ t\in[t_{0,}t_{1}]\subset[a,b],
\]
which maps initial particle positions $x_{0}$ at time $t_{0}$ to
their time-$t$ positions $x(t;t_{0},x_{0})$. 

The deformation gradient $\nabla F_{t_{0}}^{t}$ governs the infinitesimal deformations of
the phase space $U$. In particular, an infinitesimal perturbation
$\zeta_{0}$ at point $x_{0}$ and time $t_{0}$ is mapped, under
the system \eqref{eq:FlowODE} to its time-$t$ position, $\zeta_{t}=\nabla F_{t_{0}}^{t}(x_{0})\zeta_{0}$.
The squared magnitude of the evolving perturbation is governed by 

\begin{equation}
\langle\zeta_{t},\zeta_{t}\rangle=\langle\zeta_{0},C_{t_{0}}^{t}(x_{0})\zeta_{0}\rangle,\ \ \ C_{t_{0}}^{t}(x_{0})=\left[\nabla F_{t_{0}}^{t}(x_{0})\right]^{\top}\nabla F_{t_{0}}^{t}(x_{0}),\label{eq:CGdef}
\end{equation}
where $C{}_{t_{0}}^{t}$ denotes the right Cauchy--Green strain tensor \cite{TruesdellNoll2004}.
The eigenvalues $\lambda_{i}(x_{0})$ and eigenvectors $\xi_{i}(x_{0})$
of $C_{t_{0}}^{t}(x_{0})$ are defined, indexed and oriented here
through the relationship
\[
C_{t_{0}}^{t}(x_{0})\lambda_{i}(x_{0})=\lambda_{i}(x_{0})\xi_{i}(x_{0}),\qquad\left|\xi_{i}\right|=1,\ \ i=1,2;\ \ \lambda_{1}\leq\lambda_{2},\quad\xi_{2}=R\xi_{1}.
\]
For notational simplicity, we omit the dependence of $\lambda_{i}(x_{0})$
and $\xi_{i}(x_{0})$ on $t_{0}$ and $t$.

Objective coherent structures are defined as stationary curves of objective (frame-invariant)
variational principles, and can be viewed also as null-geodesics of
suitably defined Lorentzian metrics, with specific boundary conditions
\cite{LCSHallerAnnRev2015,SerraHaller2015}. These metrics are summarized
in Table \ref{SummaryTableMetricOCSs}. 

\begin{table}[h]
\begin{centering}
\begin{tabular}{l|c|c}
\hline 
\multirow{2}{*}{$\mathbf{Type\ of\ OCS}$} & \multicolumn{2}{c}{$\mathbf{Metric}:\ g(u,u)=\langle u,Au\rangle$}\tabularnewline
\cline{2-3} 
 & $\mathbf{LCS}$ & $\mathbf{OECS}$\tabularnewline
\hline 
\multirow{2}{*}{Hyperbolic \& Parabolic} & \multirow{2}{*}{$A=\tfrac{1}{2}[C_{t_{0}}^{t}R-RC_{t_{0}}^{t}]$} & \multirow{2}{*}{$A=2SR$}\tabularnewline
 &  & \tabularnewline
\hline 
\multirow{1}{*}{Elliptic } & $A_{\lambda}=\tfrac{1}{2}[C_{t_{0}}^{t}-\lambda^{2}I],\ \ \lambda\in\mathbb{R}$ & $A_{\mu}=S-\mu I,\ \ \mu\in\mathbb{R}$\tabularnewline
\hline 
\end{tabular}
\par\end{centering}

\caption{Lorentzian metrics whose null-geodesics define various coherent structures (see \cite{LCSHallerAnnRev2015,SerraHaller2015} for a review.)}
\centering{}\label{SummaryTableMetricOCSs}
\end{table}
Although eq. \eqref{eq:EulLagr_On_0LS_Pc-1-1}
can generally be applied to compute all the coherent structures listed in Table \ref{SummaryTableMetricOCSs}, here we focus
on elliptic OCSs. Elliptic OCSs are closed null-geodesics of the corresponding
Lorentzian metric families shown in Table \ref{SummaryTableMetricOCSs}.
In fluid dynamics terms, elliptic LCSs are exceptionally
coherent vortex boundaries that show no unevenness in their
tangential deformation. Similarly, elliptic OECSs are exceptionally
coherent vortex boundaries that show no infinitesimally short-term unevenness in their
tangential deformation. The parameter $\lambda$ represents the tangential
stretching experienced by an elliptic LCS over the time interval $[t_{0},t]$,
while $\mu$ denotes the tangential stretch rate along an elliptic
OECS. In the next two sections, applying Theorem \ref{sec:Thm1}
to the Lorentzian metric families $A_{\mu}$ and $A_{\lambda}$,
we derive initial value problems (ODEs and initial conditions) for the computation of Eulerian and Lagrangian vortex boundaries.

\selectlanguage{english}%

\subsection{\textcolor{black}{Elliptic OECSs }}\label{sec:EllOECSs}

Elliptic OECSs are closed null-geodesics of the
one-parameter family of Lorentzian metrics (cf. Table \ref{SummaryTableMetricOCSs})
\begin{equation*}
A_{\mu}(x,t)=S(x,t)-\mu I.
\end{equation*}
We denote by $S^{ij}(x)$
the entry at row $i$ and column $j$ of $S(x,t)$ at a fixed time
$t$, and its derivatives $\partial_{(\cdot)}S^{ij}(x)$
by $S_{(\cdot)}^{ij}(x)$. A direct application
of Theorem \ref{sec:Thm1}, leads to the the following
result.

At each time $t$ and for a given value of $\mu$, elliptic
OECSs satisfy the pointwise condition 
\begin{equation}
S^{11}(x,t)\cos^{2}\phi+S^{12}(x,t)\sin 2\phi+S^{22}(x,t)\sin^{2}\phi-\mu=0,\ \ \  (x,\phi)\in V_t,\label{eq:LevelSurfaceThmOECSs}
\end{equation}
with the set $V_t$ defined as
\begin{equation}
V_t=\{x\in U,\ \phi\in\mathbb{S}^{1}:\ \sin 2\phi[S^{22}(x,t)-S^{11}(x,t)]+2\cos 2\phi S^{12}(x,t)\neq0\}.
\label{eq:DoEOECSs}
\end{equation}
Elliptic OECSs coincide with the $x-$projection of closed orbits
of the initial value problem
\begin{equation}
\begin{aligned}
x^{\prime} & = e_{\phi},\\
\phi^{\prime} & = -\frac{\cos^{2}\phi\langle\nabla_{x}S^{11}(x,t),e_{\phi}\rangle+\sin 2\phi\langle\nabla_{x}S^{12}(x,t),e_{\phi}\rangle+\sin^{2}\phi\langle\nabla_{x}S^{22}(x,t),e_{\phi}\rangle}{\sin 2\phi[S^{22}(x,t)-S^{11}(x,t)]+2\cos 2\phi S^{12}(x,t)},\label{eq:CorollEllOECSs}
\end{aligned}
\end{equation}
\begin{equation}
(x_{0},\phi_{0})=\left\{ \left(x_{0}(\mu,0),0\right)\in V_t:\ \ S^{11}(x_0)-\mu=0\right\}.\label{eq:r0OECSs}
\end{equation}

\noindent In the case of incompressible flows $(\nabla\cdot f\equiv0)$, eq. \eqref{eq:CorollEllOECSs} simplifies to 
\begin{equation}
\begin{aligned}
x^{\prime} & = e_{\phi},\\
\phi^{\prime} & =  -\frac{[S_{x_{1}}^{11}(x,t)\cos\phi+S_{x_{2}}^{11}(x,t)\sin\phi]\cos 2\phi +[S_{x_{1}}^{12}(x,t)\cos\phi+S_{x_{2}}^{12}(x,t)\sin\phi]\sin 2\phi}{2[S^{12}(x,t)\cos 2\phi-S^{11}(x,t)\sin 2\phi ]}.\label{eq:CorollEllOECSsIncompr}
\end{aligned}
\end{equation}

\subsubsection{Elliptic OECSs: streamfunction formulation}\label{sec:CorEllOECSsHamiltonian}

In case the velocity fields is derived from a time-dependent streamfunction
$\psi(x,t)$, the ODE \eqref{eq:FlowODE} is of the form
\begin{equation}
\begin{aligned}
\dot{x}_{1} & = -\psi_{x_{2}}(x_{1},x_{2},t)\\
\dot{x}_{2} & = \psi_{x_{1}}(x_{1},x_{2},t).\label{eq:HAmiltVelFields}
\end{aligned}
\end{equation}
Denoting the partial derivative $\partial_{x_{i}}\psi(x)$ by $\psi_{i}(x),\ i\in$\{1,2\},
we reformulate our results in terms of the time-dependent streamfunction as follows.

For a velocity field generated by the time-dependent
streamfunction $\psi(x_{1},x_{2},t)$, 
at each time $t$ and for a given value of $\mu$, elliptic
OECSs satisfy the pointwise condition
\begin{equation}
\psi_{21}(x,t)\cos 2\phi+\tfrac{1}{2}[\psi_{22}(x,t)-\psi_{11}(x,t)]\sin 2\phi +\mu=0,\ \ \ (x,\phi)\in V_t,\label{eq:LevelSurfaceThm_Hamilt}
\end{equation}
within the set $V_t$ defined as 
\begin{equation}
V_t=\{x\in U,\ \phi\in\mathbb{S}^{1}:\ [\psi_{11}(x,t)-\psi_{22}(x,t)]\cos 2\phi+2\psi_{21}(x,t)\sin 2\phi\neq0\}.
\label{eq:DoEOECSsStrFun}
\end{equation}
Furthermore, elliptic OECSs coincide with the $x-$projection of closed orbits of the initial value problem
\begin{equation}
\begin{aligned}
x^{\prime} & =e_{\phi}, \\
\phi^{\prime} & = -\frac{\langle\nabla_{x}\psi_{21}(x,t),e_{\phi}\rangle\cos 2\phi +\tfrac{1}{2}\langle\nabla_{x}[\psi_{11}(x,t)-\psi_{22}(x,t)],e_{\phi}\rangle\sin 2\phi}{[\psi_{11}(x,t)-\psi_{22}(x,t)]\cos 2\phi+2\psi_{21}(x,t)\sin 2\phi},\label{eq:CorollEllOECSsHamilt}
\end{aligned}
\end{equation}

\begin{equation}
(x_{0},\phi_{0})=\left\{ \left(x_{0}(\mu,0),0\right)\in V_t:\ \ \psi_{21}(x_0,t)-\mu=0\right\}.\label{eq:r0Hamilt}
\end{equation}

\subsection{Elliptic LCSs\foreignlanguage{american}{\label{sec:CorEllLCSs}}}

For Lagrangian vortex boundaries
(elliptic LCSs), the underlying Lorentzian metric is (cf. Table \ref{SummaryTableMetricOCSs}) 
\begin{equation*}
A_{\lambda}(x)=C_{t_{0}}^{t}(x)-\lambda^{2}I.
\end{equation*}
To avoid confusion with the initial conditions of the reduced null-geodesic flow (cf. eq.\eqref{eq:r0alpha-1}), here we denote the spatial dependence of the Cauchy-Green by $x$ instead of $x_0$.  Applying Theorem \ref{sec:Thm1}, and denoting by $C^{ij}(x)$ the entry at row $i$ and column $j$ of $C_{t_{0}}^{t}(x)$ we obtain the following result.

For a fixed time interval $[t_{0},t_1]$ and for a given value of $\lambda$,
elliptic LCSs satisfies pointwise the condition
\begin{equation}
C^{11}(x)\cos^{2}\phi+C^{12}(x)\sin 2\phi +C^{22}(x)\sin^{2}\phi-\lambda^{2}=0,\ \ \  (x,\phi)\in V,\label{eq:LevelSurfaceThm_LCS}
\end{equation}
within the set $V$ defined as
\begin{equation}
V=\{x\in U,\ \phi\in\mathbb{S}^{1}:\ \sin 2\phi [C^{22}(x)-C^{11}(x)]+2\cos 2\phi C^{12}(x)\neq0\}.
\label{eq:DoEOLCSs}
\end{equation} 
Furthermore, elliptic LCSs coincide with the $x-$projection of closed orbits
of the initial value problem 
\begin{equation}
\begin{aligned}
x^{\prime} & = e_{\phi}, \\
\phi^{\prime} & = -\frac{\cos^{2}\phi\langle\nabla_{x}C^{11}(x),e_{\phi}\rangle+\sin 2\phi\langle\nabla_{x}C^{12}(x),e_{\phi}\rangle+\sin^{2}\phi\langle\nabla_{x}C^{22}(x),e_{\phi}\rangle}{\sin 2\phi [C^{22}(x)-C^{11}(x)]+2\cos 2\phi C^{12}(x)},\label{eq:CorollEllOECSs-2}
\end{aligned}
\end{equation}

\begin{equation}
(x_{0},\phi_{0})=\left\{ \left(x_{0}(\lambda,0),0\right)\in V:\ \ C^{11}(x_0)-\lambda^{2}=0\right\}.\label{eq:r0alpha-1-1-2}
\end{equation}

\section{Example: Mesoscale coherent vortices in large-scale ocean data\label{sec:OceanExample}}

We now use the results of Sections
\ref{sec:CorEllOECSsHamiltonian} and \ref{sec:CorEllLCSs}
to locate coherent vortex boundaries in a two-dimensional ocean-surface-velocity dataset derived from AVISO satellite altimetry measurements
(\url{http://www.aviso.oceanobs.com}). The domain of interest is
the Agulhas leakage in the Southern Ocean, bounded by longitudes $[3^{\circ}W,1^{\circ}E]$,
latitudes $[32^{\circ}S,24^{\circ}S]$ and the time slice we selected
correspond to $t=24\ \mathrm{November}\ 2006$. This dataset has also been used in the vortex detection studies \cite{HallerBeron-Vera2012,Wang2015a,SerraHaller2015}, which provide a benchmark for comparison with the approach developed here.

Under the geostrophic assumption, the ocean surface height measured
by satellites plays the role of a streamfunction for surface currents.
With $h$ denoting the sea surface height, the velocity field in longitude-latitude coordinates, $[\varphi,\theta]$, can be expressed as 
\[
\dot{\varphi}=-\dfrac{g}{R^{2}f_c(\theta)\cos\theta}\partial_{\theta}h(\varphi,\theta,t),\ \ \ \dot{\theta}=\dfrac{g}{R^{2}f_c(\theta)\cos\theta}\partial_{\varphi}h(\varphi,\theta,t),
\]
where $f_c(\theta):=2\Omega\sin\theta$ denotes the Coriolis parameter,
$g$ the constant of gravity, $R$ the mean radius of the earth and
$\Omega$ its mean angular velocity. The velocity field is available
at weekly intervals, with a spatial longitude-latitude resolution
of $0.25^{\circ}$. For more detail on the data, see \cite{Beron-VeraDatassh2013}. 

\selectlanguage{english}%

\subsection{\textcolor{black}{Elliptic OECSs \label{sec:OceanExampleE}}}

\begin{figure}[h]
\subfloat[\selectlanguage{english}%
\selectlanguage{american}%
]{\includegraphics[height=0.43\columnwidth]{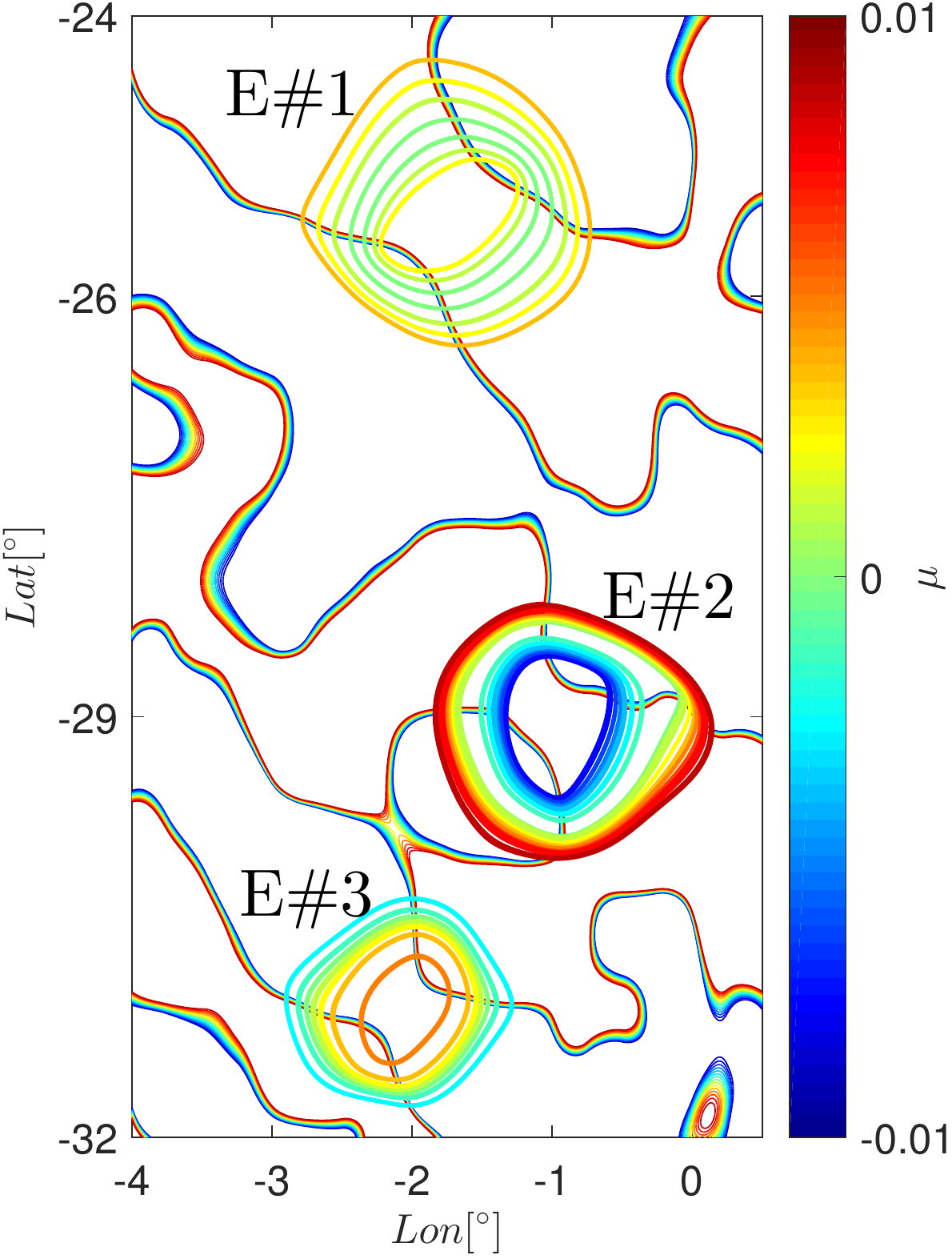}\label{fig:OECSsonIC}}\hfill{}\subfloat[\selectlanguage{english}%
\selectlanguage{american}%
]{\includegraphics[height=0.43\columnwidth]{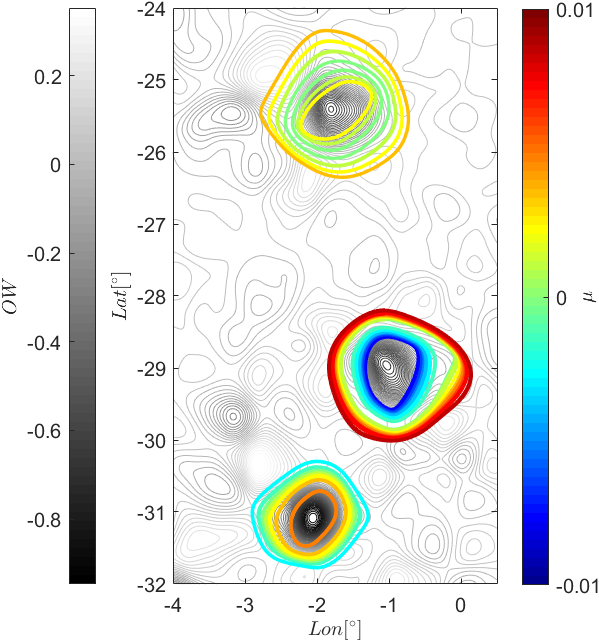}

\label{fig:OECSsonOW}}\hfill{}\subfloat[\selectlanguage{english}%
\selectlanguage{american}%
]{\includegraphics[height=0.3\textheight]{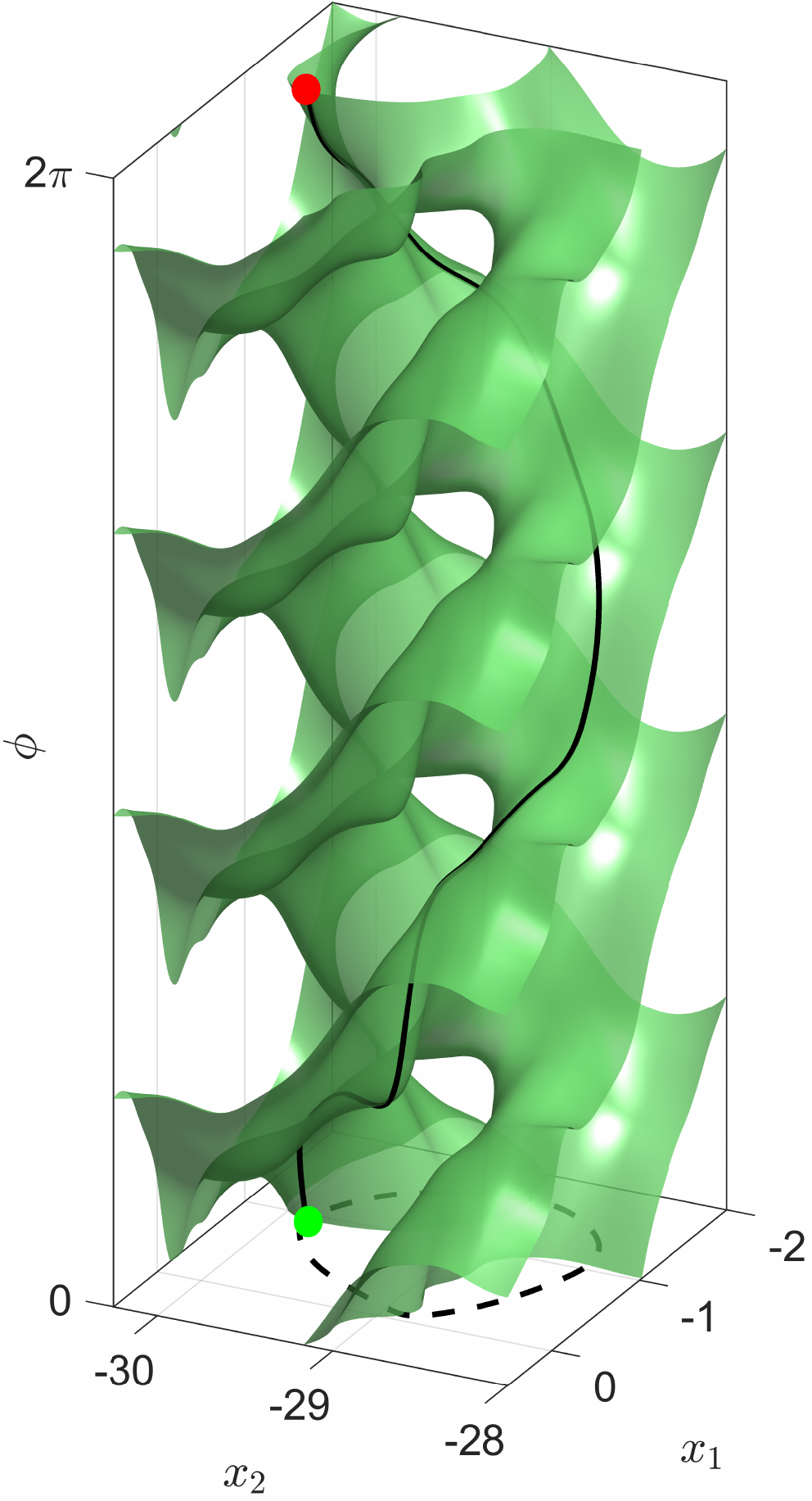}

\label{fig:r(s)onG_const}}

\caption{(a) Elliptic OECSs for different values of the stretching rate $\mu$
(in color) along with the $x-$component of the initial conditions
$x_{0}(\mu,0)$, as defined in eq. \eqref{eq:r0Hamilt}.
(b) The same elliptic OECSs of Fig. \ref{fig:OECSsonIC} on level sets of the OW parameter encoded with
the gray colormap. (c) Elliptic OECS for
$\mu=0$ corresponding to the vortical region denoted by E$\#2$ (cf. Fig \ref{fig:OECSsonIC}), in
the $\mathbb{R}^{2}\times\mathbb{S}^{1}$ space. The green surface represents the
zero set defined by eq. \eqref{eq:LevelSurfaceThm_Hamilt}, the solid
black curve is the periodic orbit of the ODE \eqref{eq:CorollEllOECSsHamilt},
and the dashed black curve is the corresponding elliptic OECSs.}

\label{fig:OECSsOcean}
\end{figure}

\selectlanguage{american}%
Applying the results in Section \foreignlanguage{english}{\textcolor{black}{\ref{sec:CorEllOECSsHamiltonian}}},
we obtain three objectively detected vortical regions in the domain
under study, each filled with families of elliptic OECSs (cf. Fig.
\foreignlanguage{english}{\textcolor{black}{\ref{fig:OECSsOcean}}}).
\foreignlanguage{english}{Figure \textcolor{black}{\ref{fig:OECSsonIC}
shows elliptic OECSs for different values of stretching rate $\mu$
(in color), along with the $x-$component of the initial conditions
$x_{0}(\mu,\phi_{0})$ for $\phi_{0}=0$ (cf. eq. \eqref{eq:r0Hamilt}
or Fig. \ref{fig:InitialCond}). }Figure \textcolor{black}{\ref{fig:OECSsonOW}
shows the same elliptic OECSs of Fig. \ref{fig:OECSsonIC} }}along
with level sets of the Okubo--Weiss (OW) parameter
\[
OW(x,t)=s_{2}^{2}(x,t)-\omega^{2}(x,t),
\]
where $\omega(x,t)$ denotes the vorticity. Spatial domains with $OW(x,t)<0$
are frequently used indicators of instantaneous ellipticity in unsteady
fluid flows \cite{Okubo1970,Weiss1991}. The OW parameter, however,
is not objective (the vorticity term will change under rotations),
and can hence generate both false positives and false negatives in vortex detection (see e.g., \cite{SerraHaller2015}).

Figure \ref{fig:r(s)onG_const}
shows the elliptic OECSs in correspondence of Eddy $\#2$ (cf. Fig \ref{fig:OECSsonIC}), for $\mu=0$,
in the $\mathbb{R}^{2}\times\mathbb{S}^{1}$ space. Specifically,
the green surface represent the zero set described by 
eq. \eqref{eq:LevelSurfaceThm_Hamilt}, the solid black curve represents
the closed integral curve of the ODE \eqref{eq:CorollEllOECSsHamilt}
with boundary points $(x_{0},0)$ (green circle) and $(x_{0},2\pi)$
(red circle). Note that the $\phi=0$ and the $\phi=2\pi$
planes are identical, but for illustration purposes, we show the $\mathbb{R}^{2}\times\mathbb{S}^{1}$
space in Cartesian coordinates rather than in toroidal coordinates.
The dashed black curve represents the corresponding elliptic OECS,
i.e., the $x-$projection of the solid black curve. The domain analyzed
in Fig. \ref{fig:r(s)onG_const} is identical to
the one used for illustration in Fig. \ref{fig:InitialCond}. 

\subsection{\textcolor{black}{Elliptic LCSs }\foreignlanguage{american}{\label{sec:OceanExampleL}}}

\begin{figure}[h]
\subfloat[\selectlanguage{english}%
\selectlanguage{american}%
]{\includegraphics[height=0.4\columnwidth]{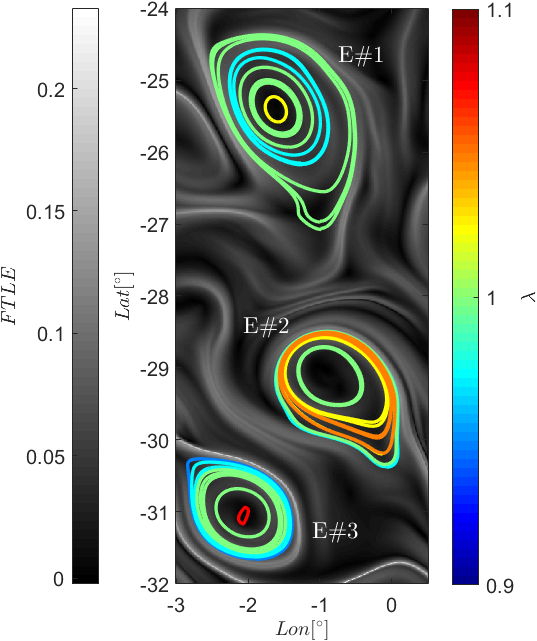}\label{fig:t0LCSs}}\hfill{}\subfloat[\selectlanguage{english}%
\selectlanguage{american}%
]{\includegraphics[height=0.4\columnwidth]{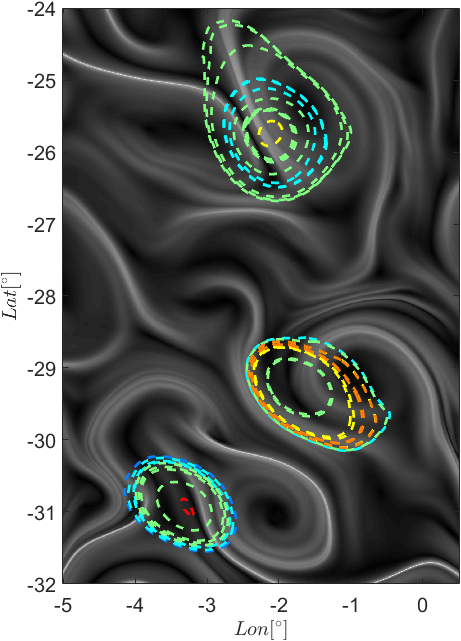}
\label{fig:t30LCSs}}\hfill{}\subfloat[\selectlanguage{english}%
\selectlanguage{american}%
]{\includegraphics[height=0.4\columnwidth]{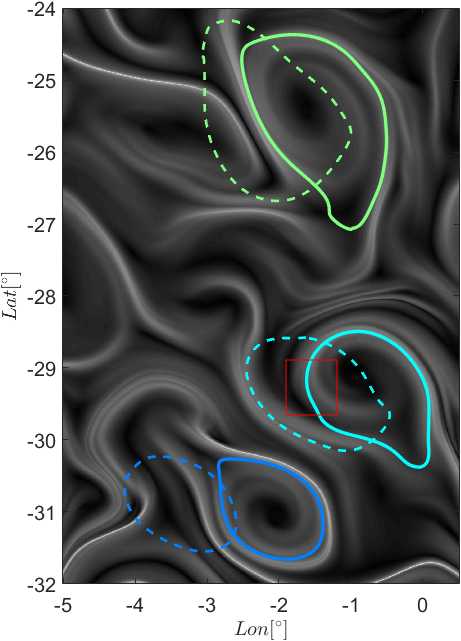}
\label{fig:t30OutLCSs}}

\caption{(a) Elliptic LCSs for different values of stretching ratio $\lambda$
(right colorbar) along with the FTLE field at $t_{0}$ (left colorbar).
(b) Advected images of the elliptic LCSs of Fig. \foreignlanguage{english}{\textcolor{black}{\ref{fig:t0LCSs}}}
at $t_{0}+30\ \text{days}$ on the FTLE field at $t_{0}$. (c) Outermost
elliptic LCSs of Fig. \foreignlanguage{english}{\textcolor{black}{\ref{fig:t0LCSs}}} (solid lines), together with their corresponding advected
images at $t_{0}+30\ \text{days}$ (dashed lines), on the FTLE field at $t_{0}$. The red square highlights a region where the local FTLE ridge crosses the elliptic LCS corresponding to Eddy $\#2$.}

\label{fig:LCSsOcean}
\end{figure}

In our Lagrangian analysis, we consider the time interval $[t_{0},t_{0}+T]$,
with\foreignlanguage{english}{\textcolor{black}{{} }}\textcolor{black}{$t_{0}=24\ \mathrm{November}\ 2006$
and}\foreignlanguage{english}{\textcolor{black}{{} $T=30\text{ days}$.}}
Applying the results in \foreignlanguage{english}{\textcolor{black}{Section
\ref{sec:CorEllLCSs},}} we obtain three objectively detected Lagrangian
coherent vortices in the domain under study, each filled with families
of elliptic LCSs (cf. Fig. \foreignlanguage{english}{\textcolor{black}{\ref{fig:LCSsOcean}}}).

\selectlanguage{english}%
Figure \ref{fig:t0LCSs} shows elliptic LCSs for
different values of the stretching ratio $\lambda$ (right colorbar),
along with the finite time Lyapunov exponent (FTLE)
field 
\[
\Lambda(x_{0},t_{0},T)=\dfrac{1}{2T}\log(\lambda_{2}(x_{0},t_{0},t_{0}+T)),
\]
encoded with the left colorbar. The FTLE measures the maximal local
separation of nearby initial conditions over the time interval $[t_{0,}t_{0}+T]$. FTLE ridges are usually used as a visual diagnostic to distinguish coherent regions from the surrounding chaotic regions. The FTLE field, however, does not give any vortex
	boundary, and can incorrectly indicate the presence of LCSs \cite{LCSHallerAnnRev2015}.
	Moreover, the extraction of FTLE ridges requires sophisticated
	post-processing algorithms (see e.g., \cite{Eberly1996}).
	This is mainly because ridges separate regions of the phase space
	with different behaviors, increasing considerably the
	sensitivity of any numerical computation in their vicinity. Examples
	include the detection of Cauchy--Green singularities, which plays
	a crucial role in direction-field-based procedure for the computation of elliptic
	LCSs (cf. Appendix \ref{sec:DirFieldApproachIC} or \cite{AutoMdetectFlorian}). Specifically, near FTLE ridges, singularities tend to artificially cluster (cf. Fig. \ref{fig:ClustCGsing}) preventing, possibly, the identification of the outermost elliptic LCSs.
	
On this note, Fig. \ref{fig:t30OutLCSs} shows that the initial position (solid line) of the outermost elliptic LCSs in correspondence of Eddies $\#1$ and $\#3$, almost overlap with nearby FTLE ridges. In contrary, the outermost elliptic LCSs in correspondence of Eddy $\#2$ crosses the local FTLE ridge (cf. red square in Fig. \ref{fig:t30OutLCSs}). The dashed lines represent the final position of the outermost elliptic LCSs. This highlights two important facts. First, elliptic LCSs computed with the present scheme are insensitive to artificial clusters of singularities, and hence identify the correct  boundary of coherent Lagrangian vortices. Second, FTLE ridges do not signal correct Lagrangian vortex boundaries. 

Figure \ref{fig:t30LCSs} shows the advected images of elliptic LCSs of Fig. \ref{fig:t0LCSs} at time $t_{0}+30\ \text{days}$, along with the FTLE field at $t_{0}$. All vortex boundaries remain perfectly coherent for a time interval equal to the extraction time $T$, as expected.

\section{Conclusions}

Recently developed variational methods offer exact
definitions for Objective Coherent Structures (OCSs) as centerpieces
of observed trajectory patterns.
OCSs can be classified into Lagrangian Coherent Structures (LCSs) \foreignlanguage{american}{\cite{LCSHallerAnnRev2015}}
and Objective Eulerian Coherent Structures
(OECSs) \cite{SerraHaller2015}, depending on the time interval over
which they shape trajectory patterns. LCSs
are intrinsically tied to a specific finite time interval over which
they are influential, while OECSs are computable
at any time instant, with their influence confined to short time
scales. Both type of OCSs can be computed as null-geodesics of suitably
defined Lorentzian metrics defined on the physical domain of the underlying
fluid. 

Prior numerical procedures for the computation of such vortex
boundaries require significant numerical effort to overcome the sensitivity of the steps involved.
Here we have derived and tested a simplified and more accurate numerical method. Our method is based on a direct solution of a reduced, three-dimensional version of the underlying ODEs for null-geodesics. Based on topological properties of simple planar closed curves, we also derive the admissible set of initial conditions for the reduced ODEs overcoming the limitation of the existing procedure, and making the detection of null-geodesic fully automated. In the supplementary material, we provide a MATLAB implementation of this method, with further explanation in Appendix \ref{sec:Algorithm}.

 We have illustrated the present method on mesoscale eddy-boundary extraction from satellite-inferred
ocean velocity data.

\section*{Supplementary material\label{sec:SuppMat}}

A MATLAB code for the computation of closed null-geodesics is available at \url{https://github.com/MattiaSerra/Closed-Null-Geodesics-2D}. Specifically, the MATLAB code computes elliptic LCSs (cf. Section \ref{sec:CorEllLCSs}). Appendix \ref{sec:Algorithm} summarizes the different steps of the main code with explicit references to the different subfunctions.

\section*{Acknowledgment}

We would like to acknowledge Alireza Hadjighasem for helpful discussions
on the development of the MATLAB code available as supplementary material.
\clearpage

\appendix
\section{Direction field approach for computing elliptic LCSs} \label{sec:DirFieldApproach}
Using the notation introduced in Section \ref{sec:setupAndNot}, we summarize here the direction-field approach for the computation of elliptic LCSs derived in \cite{BlackHoleHaller2013}. 
The initial position of elliptic LCSs coincide with limit cycles of the differential equation family
	\begin{equation}
	\frac{dx}{ds}=\eta_{\lambda}^{\pm}(x),\qquad\eta_{\lambda}^{\pm}=\sqrt{\dfrac{\lambda_{2}-\lambda^2}{\lambda_{2}-\lambda_{1}}}\xi_{1}\ \pm\ \sqrt{\dfrac{\lambda^2-\lambda_{1}}{\lambda_{2}-\lambda_{1}}}\xi_{2}.\label{eq:diffieldODE}
	\end{equation}
The direction field $\eta_\lambda ^\pm (x)$	depend explicitly on $\lambda$, and due to the lack of a well-defined orientation for eigenvector fields, it is a-priori unknown which one of the $\eta_\lambda ^\pm (x)$ fields can have limit cycles. Therefore, both direction fields ($\pm$) should be checked. Similar arguments and expressions hold for elliptic OECSs \cite{SerraHaller2015}.
\subsection{Selection of initial conditions} \label{sec:DirFieldApproachIC}
Here we summarize an automated procedure for the selection of initial conditions (or Poincaré Sections) of \eqref{eq:diffieldODE}, developed in \cite{AutoMdetectFlorian}. Such procedure is based on the location of Cauchy-Green singularities, whose identification is a highly sensitive procedure. 
\begin{figure}[h]
	\subfloat[]{\includegraphics[width=0.4\columnwidth]{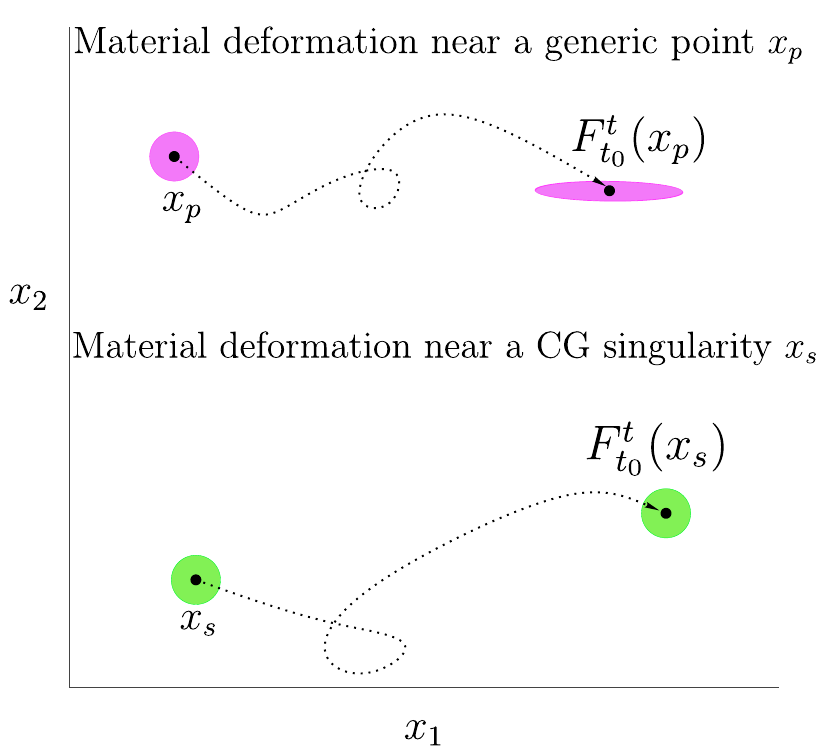}\label{fig:DefatSing}}\hfill{}
	\subfloat[]{\includegraphics[width=0.5\columnwidth]{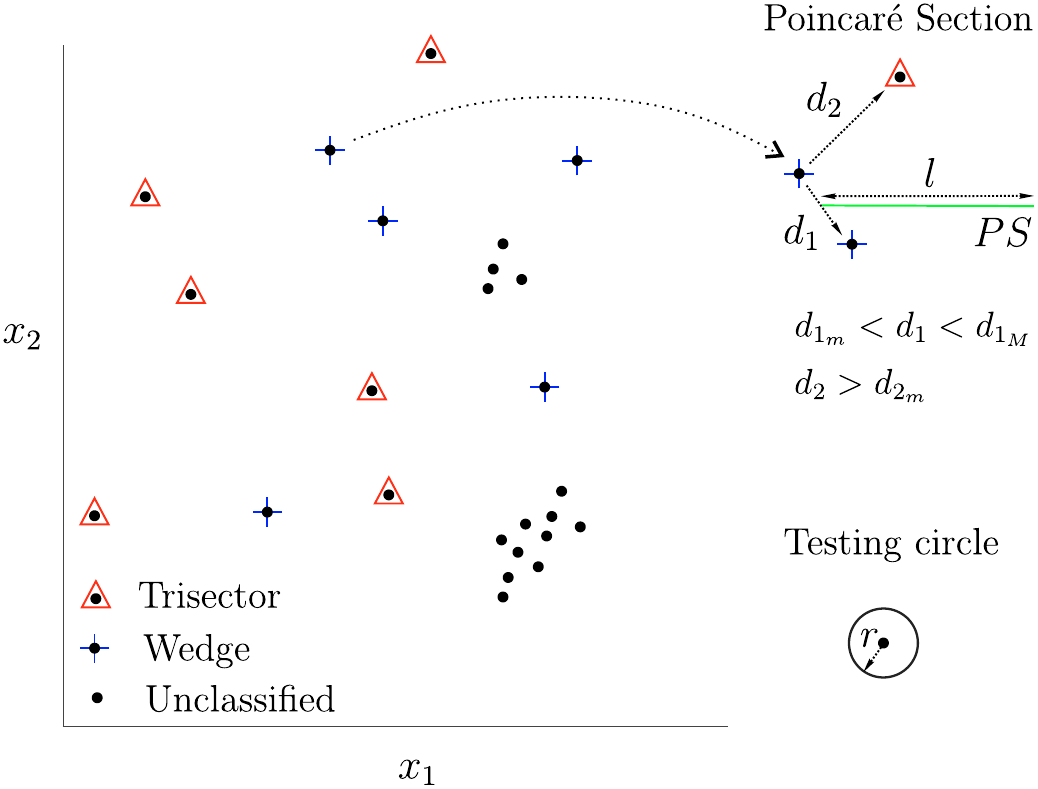}		
		\label{fig:AutomDetect}}\caption{(a) Material deformation in the neighborhood of a generic point $x_p$, and of a singularity of the Cauchy-Green tensor $x_s$, over a finite time interval $[t_0,t]$. (b) Identification of the topological type of tensor-line singularities, and definition of the Poincaré Section (PS) for the direction field integration, as in \cite{AutoMdetectFlorian}. The user-input parameters $r,l,d_{1_m},d_{1_M},d_{2_m}$ are used to locate the PS, and depend on the specific problem. Specifically, $r$ is the radius of the testing circle used to identify the singularity type, $l$ is the length of the PS, and $d_{1_m},d_{1_M},d_{2_m}$ bounds the distances of the first two closest singularities from the selected one.}
	\label{fig:DefsingAndAutomDetect}
\end{figure}

\begin{figure}[h]
	\centering
	\subfloat[]{\includegraphics[width=0.7\columnwidth]{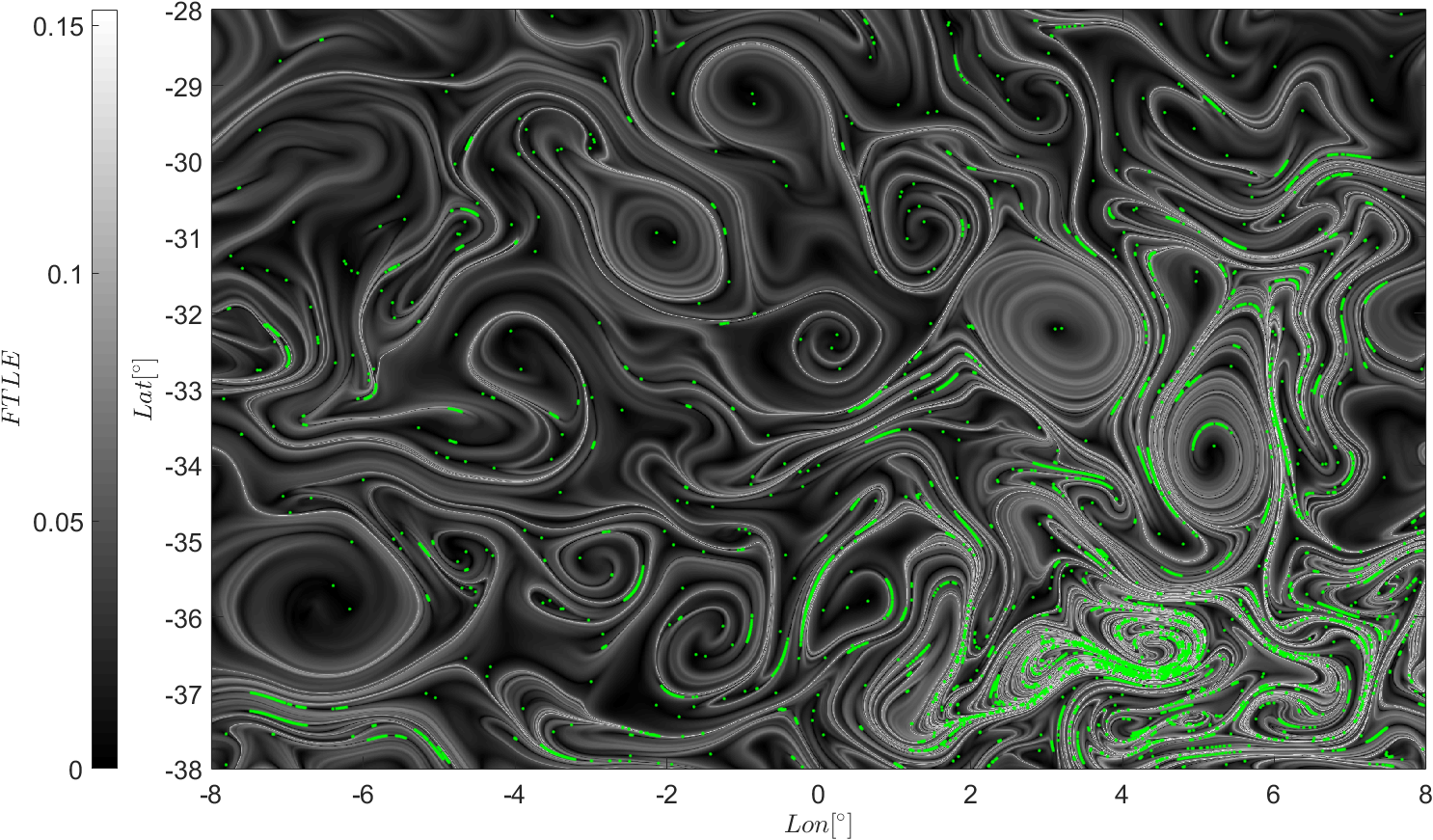}\label{fig:SingClust60d}}\\ 
	\subfloat[]{\includegraphics[width=0.7\columnwidth]{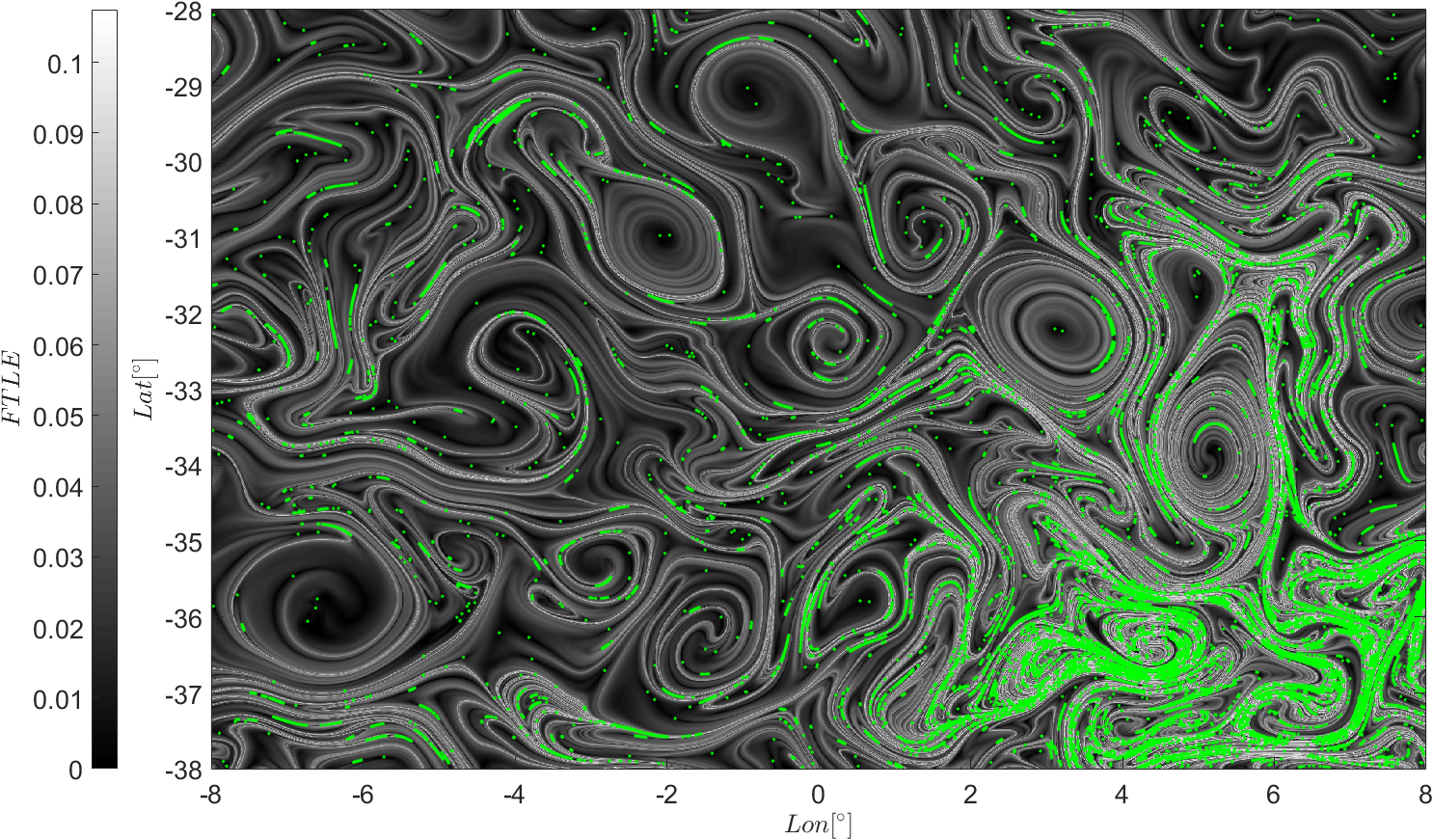}\label{fig:SingClust90d}}		
	\caption{(a) Cauchy-Green singularities (green dots) along with the FTLE field for an integration time of two months. (b) Cauchy-Green singularities (green dots) along with the FTLE field for an integration time of three months.}
	\label{fig:ClustCGsing}
\end{figure}

 Singularities of the Cauchy-Green tensor are exceptional points in the initial configuration of the fluid domain where no distinguished stretching directions exist, and hence an initially circular neighborhood around them will remain undeformed under the action of the flow map. Figure \ref{fig:DefatSing} shows the material deformation in the neighborhood of a generic point $x_p$, and of a singularity of the Cauchy-Green tensor $x_s$, over a finite time interval $[t_0,t]$.
 In a typical turbulent flow, one expects that the occurrence of these points decreases with longer time interval due to the increased mixing in the flow. The detection of tensor filed singularities, however, is a particularly sensitive process, and this sensitivity increases with longer integration times, leading to artificial clusters of singularities (cf. \textcolor{black}{Fig. \ref{fig:ClustCGsing}} or \cite{AutoMdetectFlorian}). 

In Fig. \ref{fig:AutomDetect}, we illustrate the main steps used in \cite{AutoMdetectFlorian} to locate the Poincaré Section (PS) for null-geodesics computations with the direction field approach. First, Karrash et al. identify the topological type of each singularity using a testing circle of radius $r$. When singularities are too close to each other (i.e., distance smaller than $r$), their topological type cannot be identified and they remain unclassified. Applying an index theory argument to direction fields, Karrash et al. \cite{AutoMdetectFlorian} show that each null-geodesic on ($U,g_x$) contains at least two wedge-type singularities in its interior. Relying on this necessary condition, they seek isolated wedge-pairs and set a PS of length $l$ from theirs mid points. Specifically, an isolated wedge pair exists if the distance $d_1$ between the current wedge and the closest one is such that $d_{1_m}<d_1<d_{1_M}$, and the second closest singularity is a trisector, whose distance $d_{2}>d_{2_m}$. This procedure require user-input parameters $r,l,d_{1_m},d_{1_M},d_{2_m}$, which are problem dependent. At the same time, it will also miss null-geodesics with more than one wedge pair in their interior.

Figure \ref{fig:ClustCGsing} shows the singularities of $C_{t_0}^{t_0+T}(x_0)$ (green dots), along with the corresponding FTLE field, in the flow domain bounded by longitudes $[8^{\circ}W,8^{\circ}E]$, latitudes $[38^{\circ}S,28^{\circ}S]$, with $t_0=24\ \mathrm{November}\ 2006$. Specifically, in Fig. \ref{fig:SingClust60d} the integration time $T=2$ months, while in Fig. \ref{fig:SingClust90d} $T=3$ months. This figure shows an artificial clustering of singularities with increasing integration times, which makes singularity-based methods for detecting null-geodesics non-optimal. Even if no such clustering occurs (e.g., in the detection of elliptic OECSs), a parameter-free method, as the one developed here, is preferable.  

\section{Steps for the computation of closed null-geodesics} \label{sec:Algorithm}
Algorithm \ref{algorithm1} provides a
brief summary of the main steps performed by the MATLAB code (cf. supplementary material) for the computation of closed null-geodesics. Specifically, Algorithm \ref{algorithm1} computes elliptic LCSs. We list the \textsl{MATLAB subfunctions} used to compute the formulas in Section \ref{sec:CorEllLCSs}.

\begin{algorithm}[H]
\caption{Compute elliptic LCSs (cf. Section \foreignlanguage{english}{\ref{sec:CorEllLCSs}})}
\label{algorithm1} 
\textbf{Input:} (i) Entries of the Cauchy-Green tensor
field $C^{ij}(x)$ and their spatial derivatives $C^{ij}_{x_{k}}(x),\ i,j,k\in\text{{\{}}1,2\text{{\}}}$, along with the
corresponding spatial grid vectors $xi\_g,\ i\in\text{{\{}}1,2\text{{\}}}.$
(ii) A vector $lamV$ containing the desired values of the parameter $\lambda$. 
\begin{enumerate}
\item Compute $r_{\lambda}(0)$ (cf. eq. \eqref{eq:r0alpha-1-1-2}):
\textsl{\textcolor{black}{r0$\_$lam.m}}
\item Compute $\phi^{\prime}(x_{1},x_{2},\phi)$ (cf. eq. \eqref{eq:CorollEllOECSs-2}):
\textsl{\textcolor{black}{Phi\_prime.m}}
\item Find closed null-geodesics (cf. eq. \eqref{eq:CorollEllOECSs-2}-\eqref{eq:r0alpha-1-1-2}):
\textsl{\textcolor{black}{FindClosedNullGeod.m}}
\item Find outermost closed null-geodesics: \textsl{FindOutermost.m}\textsl{\textcolor{black}{{} }}
\end{enumerate}
\textbf{Output:} Elliptic LCSs corresponding to the different values
of $\lambda$. 
\end{algorithm}
Algorithm \ref{algorithm1} is
general and can be used to compute elliptic OECSs (cf. Section \ref{sec:EllOECSs})
or any general closed null-geodesics as defined in Theorem \ref{sec:Thm1},
where $C^{ij}(x)=A^{ij}(x),\ i,j\in\text{{\{}}1,2\text{{\}}}$. Steps
2 and 3 of Algorithm \ref{algorithm1} can
be used to compute general non closed null-geodesics of Lorentzian
metrics of the form $A_{\alpha}=A(x)-\alpha I$, from given initial
conditions.

\section{Domain of existence of the reduced geodesic flow \foreignlanguage{american}{\label{sec:DomExistGeodFlow}}}

Equation \textcolor{black}{\eqref{eq:EulLagr_On_0LS_Pc-1}} only admits
non-degenerate $r(s)=(x(s),\phi(s))^{\top}$ solutions in the set
$V\subset U\times\mathbb{\mathbb{S}}^{1}:\ B(x,\phi):=\left\langle e_{\phi},R^{\top}A(x)e_{\phi}\right\rangle \neq0.$
Note that in the set $\overline{V}$ where $B(x,\phi)=0$, equation
$g_{x}(e_\phi,e_\phi)=0$ (cf. eq. \textcolor{black}{\eqref{eq:Noether-2-1}})
does not define locally a \textcolor{black}{2-dimensional manifold}
parametrized by $(x,\phi(x))$. In fact, by the implicit function
theorem, $\phi(x)$ is defined only if $g_{x}(e_\phi,e_\phi)=0$
admits a solution and $\partial_{\phi}g_{x}(e_\phi,e_\phi)\ne0$, where
\begin{equation*}
\partial_{\phi}g_{x}(e_\phi,e_\phi)= \frac{1}{2}\left\langle e_{\phi},R^{\top}A(x)e_{\phi}\right\rangle +\frac{1}{2}\left\langle e_{\phi},A(x)Re_{\phi}\right\rangle =\left\langle e_{\phi},\text{{Sym}}(R^{\top}A(x))e_{\phi}\right\rangle =\left\langle e_{\phi},R^{\top}A(x)e_{\phi}\right\rangle .
\end{equation*}
Therefore, the set $\overline{V}$ is the union of points that satisfy
at least one of the two following conditions
\[
\begin{cases}
B(\cdot,\phi)=0\iff e_{\phi}\equiv\zeta_{i}(\cdot),\ \ \ A(\cdot)\zeta_{i}(\cdot)=\alpha_{i}(\cdot)\zeta_{i}(\cdot),\ \alpha_{i}(x)\in\mathbb{R},\  i=\text{{\{}}1,2\text{{\}}},\\
B(\phi,\cdot)=0\iff A(\cdot)\ \text{{is\ degenerate}}.
\end{cases}
\]
Equivalently,
\begin{equation*}
V = \left\{ \left(x,\phi\right)\in U\times\mathbb{\mathbb{S}}^{1}:\ A(x)e_{\phi}\nparallel e_{\phi},\ A(x)\neq \mathbf{0}\right\}.
\end{equation*}
Geometrically this means that in $\overline{V}$, there cannot be
a transverse zero of the function $g_{x}(e_\phi,e_\phi)$. Specifically,
when the first condition holds, such zero in non transverse at $x$
only for the directions $\phi(x)$ aligned with the eigenvectors of
$A(x)$. When the second condition holds, there cannot be any transverse
zero at $x$ for all $\phi$, since $A(x)$ is degenerate and no distinguished
directions exist.

\section{Hamiltonian reduction of the geodesic flow\foreignlanguage{american}{\label{sec:HamiltReduction}}}

Here we use the Hamiltonian formalism to derive a reduced geodesic
flow which is equivalent to the one derived in Sections \ref{sec:FormulationPbl}-\ref{sec:FormulationPbln2}.
With the generalized momentum $p$ defined as 
\begin{equation}
p=\frac{\partial L}{\partial x^{\prime}}=A(x)x^{\prime},\label{eq:GenerMomentum}
\end{equation}
the parametrization $x(s)$ of a geodesic $\gamma$ of the metric
$g_x(u,u)=\frac{1}{2}\left\langle u,A(x)u\right\rangle $ satisfies
the first-order system of differential equations
\begin{equation}
\begin{aligned}
x^{\prime} & = A^{-1}(x)p,\\
p^{\prime} & = -\frac{1}{2}\nabla_{x}\left\langle p,A^{-1}(x)p\right\rangle ,\label{eq:HGEL-6}
\end{aligned}
\end{equation}
which is a canonical Hamiltonian system with Hamiltonian 
\begin{equation}
H(x,p)=\frac{1}{2}\left\langle p,A^{-1}(x)p\right\rangle =L(x,x^{\prime}).\label{eq: Hamiltonian-1}
\end{equation}
This Hamiltonian is constant along all geodesics of the metric $g_x$.
In particular, if $g_x$ is Lorentzian, then null-geodesics of $g_x$
lie in the zero level surface of $H(x,p)$. As in Section \ref{sec:FormulationPbln2},
we derive a reduced form
of the Hamiltonian flow \eqref{eq:HGEL-6}, which is often referred
to as the co-geodesic flow \cite{klingenberg1995riemannian}.

\subsection{Hamiltonian reduction to a three-dimensional geodesic flow}

We introduce polar coordinates in the $p$ direction by letting
\begin{equation*}
p=\rho e_{\phi},\qquad\rho\in\mathbb{R}^{+},\qquad e_{\phi}=\left(\cos\phi,\sin\phi\right)^{\top},\qquad\phi\in\mathbb{S}^{1}.
\end{equation*}
We then rewrite eq. \eqref{eq:HGEL-6} as 
\begin{equation}
\begin{aligned}
x^{\prime} & = \rho A^{-1}(x)e_{\phi}, \\
\rho^{\prime}e_{\phi}+\phi^{\prime}\rho Re_{\phi} & = -\rho^{2}\frac{1}{2}\nabla_{x}\left\langle e_{\phi},A^{-1}(x)e_{\phi}\right\rangle ,\label{eq:HGEL-3-3}
\end{aligned}
\end{equation}
that, together with the rescaling \eqref{eq:rescale-1-2-1}, gives
\begin{equation}
\begin{aligned}
\frac{dx}{d\bar{s}} & = A^{-1}(x)e_{\phi},\\
\frac{d\rho}{d\bar{s}}\rho e_{\phi}+\frac{d\phi}{d\bar{s}}\rho^{2}Re & = -\rho^{2}\frac{1}{2}\nabla_{x}\left\langle e_{\phi},A^{-1}(x)e_{\phi}\right\rangle ,\label{eq:HGEL-4-3}
\end{aligned}
\end{equation}
or, equivalently, 
\begin{equation}
\begin{aligned}
\frac{dx}{d\bar{s}} & = A^{-1}(x)e_{\phi},\\
\frac{d\phi}{d\bar{s}} & = -\frac{1}{2}\langle \nabla_{x}\left\langle e_{\phi},A^{-1}(x)e_{\phi}\right\rangle, Re_{\phi}\rangle,\label{eq:HGEL-5-3}\\
\frac{d\rho}{d\bar{s}} & = -\rho\frac{1}{2}\langle\nabla_{x}\left\langle e_{\phi},A^{-1}(x)e_{\phi}\right\rangle,e_{\phi}\rangle.
\end{aligned}
\end{equation}
For any $\rho>0$, system (\ref{eq:HGEL-5-3}) has a three-dimensional
reduced flow
\begin{equation}
\begin{aligned}
\frac{dx}{d\bar{s}} & = A^{-1}(x)e_{\phi},\\
\frac{d\phi}{d\bar{s}} & = -\frac{1}{2}\langle \nabla_{x}\left\langle e_{\phi},A^{-1}(x)e_{\phi}\right\rangle,Re_{\phi}\rangle.\label{eq:reduced-3}
\end{aligned}
\end{equation}
Therefore, any solution of (\ref{eq:HGEL-6}) with $\rho>0$ admits
a projected flow of the form (\ref{eq:reduced-3}). This is due to
the existence of a global invariant foliation in (\ref{eq:HGEL-5-3})
that renders the $\left(x,\phi\right)$ coordinates of solutions independent
of the evolution of their $\rho$ coordinate. Closed orbits of (\ref{eq:reduced-3})
are, therefore, closed geodesics on $(U,g_x)$, even though they may
not be closed orbits of the full (\ref{eq:HGEL-5-3}). Note that the
$\phi$ component in eq. \eqref{eq:reduced-3} is the polar angle of
the generalized momentum (cf. eq. \eqref{eq:GenerMomentum}), which
is different than the $\phi$ in eq. \eqref{eq:EulLagr_On_0LS_Pc-1}.
Equation (\ref{eq:reduced-3}) does not appear to be available in
the literature. The use of the energy (as opposed to the momentum
$p$) as a coordinate appears in \cite{DelshamsValdes1999} in the
context of perturbations of closed geodesics by time-periodic potentials.
The reduced flow \eqref{eq:reduced-3} in the $(x,p)$ coordinates
is equivalent to the reduced flow \eqref{eq:EulLagr_On_0LS_Pc-1}
in the $(x,v)$ coordinates.

Geodesics can also be viewed as trajectories of \eqref{eq:HGEL-6}
contained in a constant level surface of the Hamiltonian $H(x,p).$
Null-geodesics, in particular, are contained in the level surface
\[
E_{0}=\left\{ \left(x,p\right)\in U\times\mathbb{R}^{2}:\,\,H(x,p)=0\right\} ,
\]
which in polar coordinates, for any $\rho>0$, can be rewritten as 

\[
E_{0}=\left\{ \left(x,\phi\right)\in U\times\mathbb{S}^{1}:\,\,H(x,\phi)=\frac{1}{2}\left\langle e_\phi,A^{-1}(x)e_\phi\right\rangle =0\right\} .
\]
Finally, one should select the initial conditions for the ODE \eqref{eq:reduced-3} as
$x(0)=x_{0}$ and $\phi(0)=\phi_{0}$ on $E_{0}$ to satisfy
$\left\langle e_{\phi_{0}},A^{-1}(x_{0})e_{\phi_{0}}\right\rangle =0.$

\bibliographystyle{plain}
\bibliography{ReferenceList2}

\end{document}